\documentclass{article}

\usepackage{amsmath, amssymb, amsthm}
\usepackage{graphicx}
\usepackage{enumerate}
\usepackage{tikz}
\usetikzlibrary{patterns}

\newcommand{\R}[1]{\mathbb{R}^{#1}}
\newcommand{\prf}[1][]{\noindent {\bf Proof{#1}:} }

\DeclareMathOperator{\Var}{\mathsf{Var}}
\DeclareMathOperator{\bd}{bd}
\DeclareMathOperator{\cl}{cl}

\DeclareMathOperator{\sign}{sign}
\DeclareMathOperator{\conv}{conv}

\theoremstyle{definition}
\newtheorem{defi}{Def.}
\newtheorem{exa}[defi]{Example}
\theoremstyle{plain}
\newtheorem{lemma}[defi]{Lemma}
\newtheorem{theorem}[defi]{Theorem}
\newtheorem{kor}[defi]{Corollary}
\theoremstyle{remark}
\newtheorem{rem}[defi]{Remark}

\title{Variances of surface area estimators based on pixel configuration counts}
\author{J\"urgen Kampf}

\begin{document}
\maketitle

\begin{abstract}
The surface area of a set which is only observed as a binary pixel image is often estimated by a weighted sum of pixel configurations counts. In this paper we examine these estimators in a design based setting -- we assume that the observed set is shifted uniformly randomly. Bounds for the difference between the essential supremum and the essential infimum of such an estimator are derived, which imply that the variance is in $O(t^2)$ as the lattice distance $t$ tends to zero. In particular, it is asymptotically neglectable compared to the bias. A simulation study shows that the theoretically derived convergence order is optimal in general, but further improvements are possible in special cases. 
\end{abstract}



\section{Introduction}
\begin{figure}%
\begin{center}
\begin{minipage}{2cm}
\vspace{-12cm}
\includegraphics[bb=382 277 432 323, width=1cm, clip=TRUE]{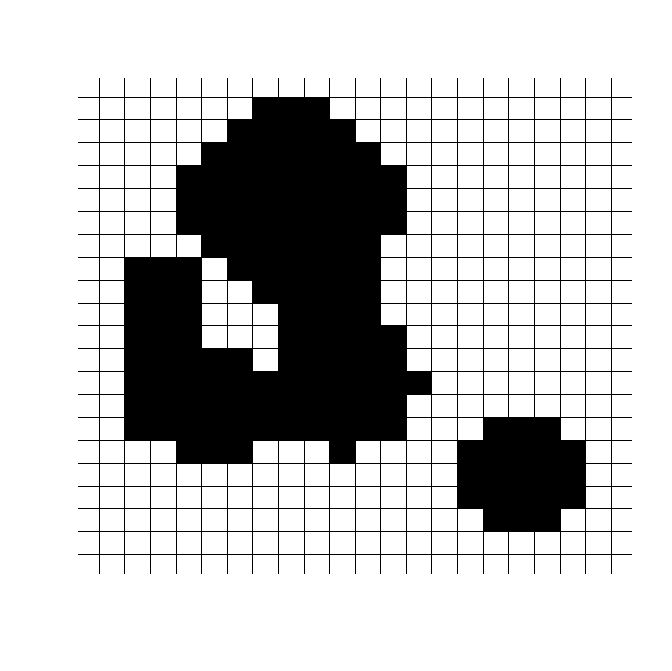}
\end{minipage}\includegraphics[width=7cm]{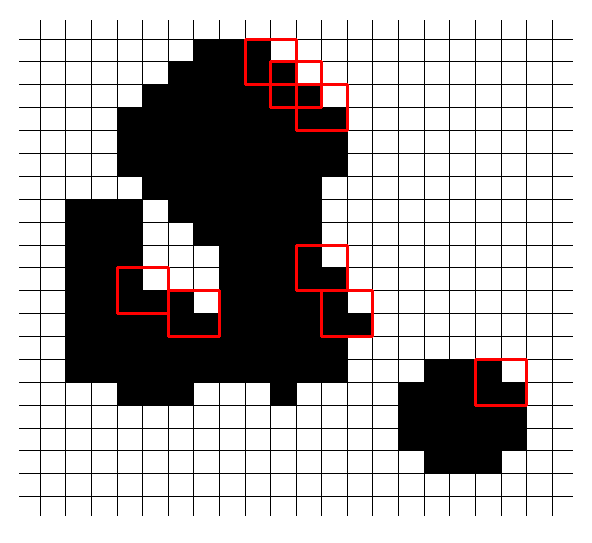}
\end{center}

\caption{Pixel configuration counts}%
\label{F:pixel_configuration_count}%
The $2\times 2$-pixel configuration on the left occurs $8$ times in the image on the right. Thus the pixel configuration count of the pixel configuration on the left in the image on the right is $8$. 

\end{figure}

There are several competing algorithms for the computation of the surface area of a set which is only observed through a pixel image, e.g.\ \cite{CFTT03, KlSu01, Li05, LiNy02, SON06}; see \cite[Sec.\ 12.2 and 12.5]{KlRo04} for an overview. A computationally fast and easy to implement approach is taken by so-called \emph{local} algorithms, \cite{Li05, LiNy02, SON06} in the above list. The idea behind these algorithms is the following: In a $d$-dimensional image a pattern of side length $n$ is called \emph{$n\times \dots \times n$-pixel configuration} ($d$ factors).  Mathematically it is modeled as a disjoint partition $(B,W)$ of $\{0,\dots, n-1\}^d$ into two disjoint subsets, where $B$ represents the set of \underline{b}lack pixels and $W$ represents the set of \underline{w}hite pixels. Since the set $\{0,\dots, n-1\}^d$ consists of $n^d$ points and each point is either colored black or colored white, there are $2^{(n^d)}$ pixel configurations. We enumerate them as $(B_j,W_j), \, j=1,\dots,2^{(n^d)}$. A pixel image of lattice distance $t>0$ can be represented by the set $A\subseteq t\mathbb{Z}^d$ of its black pixels, where $tS$ for $t>0$ and $S\subseteq \mathbb{R}^d$ is the homothetic image $\{tx\mid x\in S\}$ of $S$ with scaling factor $t$ and scaling center at the origin. While the observation window is usually bounded in applications, our results hold only if the observed set $K$ is completely contained in the observation window and thus we assume that the observation window is $\mathbb{R}^d$. Now the $j$-th \emph{pixel configuration count} at lattice distance $t$ of the image $A$ is defined as
\[ N_{t,j}(A) = \big\{ v\in\mathbb{Z}^d \mid (A-tv) \cap \{0,\dots, (n-1)t\}^d = tB_j\}. \]
It represents the number of occurrences of the $j$-th pixel configuration in the image $A$, cf.\ Figure \ref{F:pixel_configuration_count}. A local algorithm now approximates the surface area $S(K)$ of a set $K\subseteq\mathbb{R}^d$ by a weighted sum
\begin{equation} \sum_{j=1}^{2^{(n^d)}} t^{d-1} w_{j} N_{t,j}(A) \label{e:form_local} \end{equation}
of pixel configuration counts, where $A=A_t(K)$ is a pixel image of $K$ and $w_{j}$ are constants chosen in advance, called \emph{weights}. The factor $t^{d-1}$ compensates for the fact that $N_{t,j}(A)$ increases of order $(1/t)^{d-1}$ as $t\to 0$ for those pixel configurations $(B_j,W_j)$ which are ``responsible'' for the surface area (see \cite[Corollary 3.2]{Sv14} for ways how to make this precise). The two pixel configurations $(B_1, W_1)=(\emptyset, \{0,\dots, n-1\}^d)$ and $(B_{2^{(n^d)}}, W_{2^{(n^d)}})=(\{0,\dots, n-1\}^d, \emptyset)$ consisting only of white pixels resp.\ consisting only of black pixels lie typically outside resp.\ inside the image set $A$ and not on its boundary. Hence these counts do not provide information about the surface area. Thus one should put
\begin{equation} w_1=w_{2^{(n^d)}}=0, \label{e:weights_zero} \end{equation}
which will be assumed in this paper. 

For the theoretical investigation of such algorithms we assume that the set $K\subseteq\mathbb{R}^d$ is randomly shifted, i.e.\ the random set $K+tU:=\{x+tU \mid x\in K\}$ for a random vector $U$ is considered. 
A natural choice for the distribution of $U$ is the uniform distribution on $[0,1)^d$, but the results of this paper will hold of any distribution of $U$. As discretization model the \emph{Gauss discretization} is used, i.e.\ a pixel is colored black if it lies in the shifted set $K+tU$ and is colored white otherwise. Thus the model for the image of $K$ is
\[A=(K+tU) \cap  t\mathbb{Z}^d.\]
Under these assumptions no asymptotically unbiased estimator for the surface area exists. More specifically, it is shown in \cite{KiZi10} that any local estimator of the surface area in $\mathbb{R}^3$ attains relative asymptotic biases of up to $4\%$ for certain test sets, when $U$ is uniformly distributed on $[0,1)^d$, and explicit weights for which this lower bound is achieved are given. However, the bias is only one component of the error. There are a number of papers investigating the other component, namely the variance, for other estimators. Hahn and Sandau \cite{HS89} as well as Jan\'a\u{c}ek and Kub\'inov\'a \cite{JK10} investigate estimators which are suitable when the picture is analogue or when limited computational capacity requires an artificial coarsening of the image. Svane \cite{Sv15} investigates the variance of local algorithms for gray-scale images based on single pixels, i.e.\ $n=1$. There is no estimator for binary digital images for which the variance has been investigated so far. The objective of the present paper is to examine the variance of local estimators for binary digital images. 

While we use a similar setup as \cite{KiZi10}, we need slightly more strict regularity assumptions on the set $K$.  In $\mathbb{R}^2$ we assume: 
\begin{enumerate}[{(R}1)]
\item The boundary $\bd K$ of $K$ is piecewise the graph of a convex or concave function with either of the two coordinate axis being the domain, i.e.\
\[ \bd K = \bigcup_{k=1}^m F_k, \]
where for each $k\in\{1,\dots, m\}$ we have either
\[ F_k = \{ (x,f_k(x)) \mid x\in D_k \} \quad \mbox{or} \quad F_k = \{ (f_k(x),x) \mid x\in D_k \}, \]
where $D_k$ is a compact interval and $f_k: D_k \to \mathbb{R}$ is a continuous function which is convex or concave, and where $F_{k_1} \cap F_{k_2}$ for $k_1\ne k_2$ contains only points that are both endpoints of $F_{k_1}$ and of $F_{k_2}$. At the intersection point of two sets $F_{k_1}$ and $F_{k_2}$ they form an angle of strictly positive width (while an angle of $\pi$ is allowed). \label{item:R_R2}
\end{enumerate}
An even more strict assumption is needed in $\mathbb{R}^d, d>2$. 
We require:
\begin{enumerate}[{(R}1)]
\setcounter{enumi}{1}
\item The set $K$ is of the form 
\[ K = \cl\big( \bigcup_{k=1}^{m'} L_k \setminus \bigcup_{k=m'+1}^m L_k \big), \]
where $\cl S$ denotes the closure of a set $S\subseteq\mathbb{R}^d$ and where the sets $L_k$ are either convex polytopes with interior points or compact convex sets with interior points for which $\bd L_k$ is a $C^2$-manifold with nowhere vanishing Gauss-Kronecker curvature. In intersection points $z\in \bd L_{k_1} \cap \bd L_{k_2}$ the bodies $L_{k_1}$ and $L_{k_2}$ do not have a common exterior normal vector. Geometrically this means that $\bd L_{k_1}$ and $\bd L_{k_2}$ intersect nowhere under an angle of zero.  \label{item:R_Rd}
\end{enumerate}

Under the above assumptions we can show our main result. 

\begin{theorem}\label{T:main}
Let $K\subseteq \R{d}$ be a set fulfilling Assumption (R\ref{item:R_R2}) if $d=2$ or Assumption (R\ref{item:R_Rd}) if $d>2$ and let $\hat S_t$ be an estimator of the form \eqref{e:form_local} fulfilling \eqref{e:weights_zero}. Then there is a constant $s$ such that
\begin{equation}
\sup\big\{ \hat S_t(K+tv) \mid v\in[0,1)^d \big\} - \inf\big\{ \hat S_{t}(K+tv) \mid v\in[0,1)^d \big\} \le s \cdot t, \quad t>0.
\label{e:main}
\end{equation}
Thus,
\[ \Var\big(\hat S_t(K+tU) \big)\le s^2t^2/4, \quad t>0, \]
where $U$ is a random vector. 
\end{theorem}

In the first step of the proof of Theorem \ref{T:main} we show that Assumption (R\ref{item:R_R2}) or Assumption (R\ref{item:R_Rd}) implies that the boundary of $K$ can be decomposed into certain sets $M_\kappa$ to be defined below such that the intersections $M_\kappa \cap M_{\kappa'}$ are small for $\kappa\ne\kappa'$ in a certain sense. In the second step we derive upper and lower bounds for certain sums of pixel configuration counts. Since it will be possible to reconstruct the individual pixels configuration counts from these sums, the bounds derived in the second step imply the assertion of Theorem \ref{T:main}. The details are given in Section \ref{sec:proofs}. 

In Section \ref{sec:counter} we show by an example that the assertions of Theorem \ref{T:main} do not need to hold for a set $K\subseteq\mathbb{R}^d$ that is the union of two convex sets which intersect under an angle of zero. Moreover we show that an essential lemma (Lemma \ref{T:N'bounded}) of our proof fails to hold for general compact and convex sets $K\subseteq\mathbb{R}^d$. Thus the method of our proof breaks down completely without the assumption that the sets $L_1, \dots, L_m$ from (R\ref{item:R_Rd}) are either polytopes or sufficiently smooth. It is unclear, whether the assertion of Theorem \ref{T:main} still holds in this more general situation. A simulation study (Section \ref{sec:sim}) shows that the order derived in Theorem \ref{T:main} is optimal for the cube and thus is optimal in general, whereas a better bound can be achieved for the ball. In the simulation part we will also examine the integral of mean curvature. In Section \ref{sec:discuss} we discuss our results, we compare them to the results Svane \cite{Sv15} obtained for gray-scale images and we mention some open questions.

\section{The proof} \label{sec:proofs}

In this section we prove Theorem \ref{T:main}. We start by introducing some notation and in particular defining the sets $M_\kappa$ mentioned in the introduction.  Then we show that in dimension $d=2$ Assumption (R\ref{item:R_R2}) implies the existence of an appropriate boundary decomposition of $K$, followed by a proof that in dimension $d>2$ such a decomposition is implied by (R\ref{item:R_Rd}). After this, we prove that this boundary decomposition implies certain upper and lower bounds on the pixel configuration counts.  Finally we show how these bounds imply Theorem \ref{T:main}.    

\subsection{Notation} \label{ss:notation}

We assume $n$ and $d$ to be fixed and hence we will suppress dependence on $n$ and $d$ in the notation. We fix an enumeration $x_1, \dots, x_{n^d}$ of the points in $\{0,\dots, n-1\}^d$ and consider for every permutation $p\in S(n^d)$ the set
\[ \tilde G_p:= \{ u\in S^{d-1} \mid \langle x_{p(1)}, u \rangle< \langle x_{p(2)}, u \rangle < \dots < \langle x_{p(n^d)}, u \rangle \}, \]
where $S^{d-1}:=\{ u\in\mathbb{R}^d \mid \|u\| = 1 \}$, and $G_p:=\cl \tilde G_p$. Notice that 
\[ G_p:= \{ u\in S^{d-1} \mid \langle x_{p(1)}, u \rangle\le \langle x_{p(2)}, u \rangle \le \dots \le \langle x_{p(n^d)}, u \rangle \}, \] 
unless $\tilde G_p$ is empty.
If $d=2$ and $n=2$ the non-empty sets $G_p$ are the eight arcs which essentially (up to permuting the indices or changing the sign of the entries) look like
\[ G_1 = \{ (u_1, u_2)\in S^1 \mid  0 \le u_1, \ 0 \le u_2,\ u_1\le u_2\}. \]
If $n>2$ then there are more (and thus smaller) arcs.

\begin{figure}%
\begin{center}
\includegraphics[width=7cm]{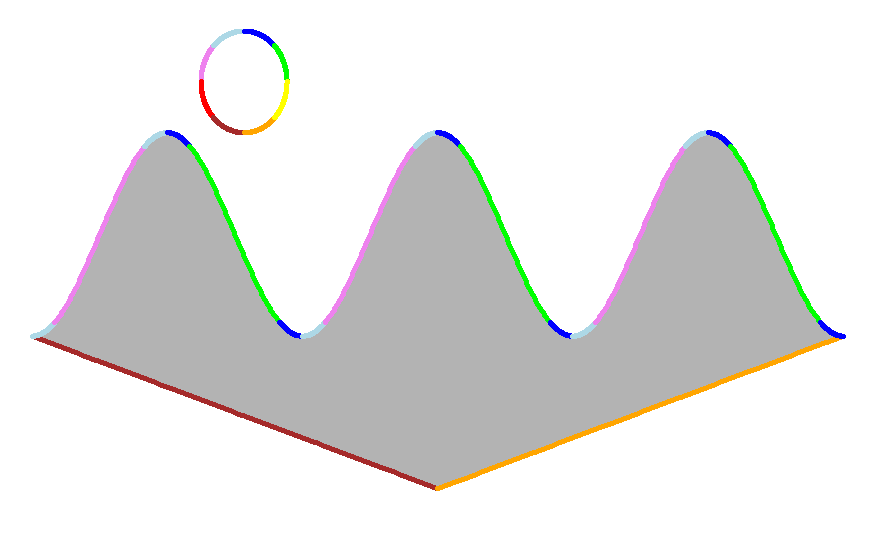}%
\end{center}
For $d=2$ and $n=2$ the unit circle $S^1$ is decomposed into eight arcs $G_p$. The boundary of each set $K\subseteq\mathbb{R}^2$ fulfilling (R\ref{item:R_R2}) - like the one printed gray here - is decomposed into eight corresponding sets $H_p$ (in this examples two of the sets $H_p$ are empty).  
\caption{The sets $G_p$ and $H_p$}%
\label{F:GpHp}%
\end{figure}

If $d=3$ and $n=2$ there are 48 sets which are isometric to
\[ G_1 = \{ (u_1, u_2, u_3) \in S^2 \mid 0\le u_1, \ 0\le u_2,\ 0\le u_3,\ u_1\le u_2, \ u_1+ u_2 \le u_3\}, \]
and 48 sets that are isometric to
\[ G_2 = \{ (u_1, u_2, u_3) \in S^2 \mid 0\le u_1, \ 0\le u_2,\ 0\le u_3, \ u_1 \le u_2\le  u_3, \ u_3 \le u_1+ u_2\}. \]
For each $p\in S(n^d)$ let $H_p$ denote the closure of the set of all boundary points of $K$ that have at least one exterior normal vector in $G_p$. Taking the closure is necessary in order to ensure $\bd K = \cup_{p\in S(n^d)} H_p$ because there are boundary points of $K$ in which there is no exterior normal vector.  An illustrative example is given in Figure \ref{F:GpHp}. 

Consider a further decomposition
\[ \bd K = \bigcup_{\kappa=1}^{\mu} M_\kappa^+, \]
where each set $M_\kappa^+$ is an intersection $F_k\cap H_p$ if $K$ fulfills (R\ref{item:R_R2}) and an intersection $(\bd L_k) \cap H_p$ if $K$ fulfills (R\ref{item:R_Rd}). In order to ensure that these sets intersect not ``not much'', put $M_{\kappa} := \cl \big(M_{\kappa}^+ \setminus \bigcup_{\lambda=1}^{\kappa-1} M_{\lambda}^+\big)$, $\kappa=1,\dots, \mu$. Let $p(\kappa), \kappa=1,\dots, \mu,$ be the element of $S(n^d)$ with $M_\kappa\subseteq H_p$ - choose an arbitrary one if it is not unique. 

A cell is a set of the form $C=[l_1, l_1+n-1] \times \dots \times [l_d, l_d+n-1]$ for some point $(l_1, \dots, l_d)\in\mathbb{Z}^d$. For $\kappa\in\{1,\dots, \mu\}$ let $\mathcal{C}^-_{\kappa,j}$ denote the system of cells $C$, such that $(\bd K) \cap C \subseteq M_\kappa$ and  $K \cap C \cap\mathbb{Z}^d = B_j+l$. Let $\mathcal{C}'_\kappa$ denote the system of all cells $C$ intersecting both $M_\kappa$ and another boundary component $M_\lambda$, $\lambda\ne \kappa$, i.e.
\begin{align*}
 \mathcal{C}'_\kappa = \{ C= [l_1, l_1+n-1] \times \dots \times [l_d, l_d+n-1] \mid& (l_1, \dots, l_d)\in\mathbb{Z}^d, \ C\cap M_\kappa \ne \emptyset, \\ 
&C\cap M_\lambda \ne \emptyset \mbox{ for some $\lambda\in\{1,\dots, \mu\}\setminus\{\kappa\}$ } \} . 
\end{align*}

Put $N^-_{\kappa,j}(K):= \# \mathcal{C}^-_{\kappa,j}$, $N_\kappa'(K):= \# \mathcal{C}'_\kappa$ and $N'(K):=\sum_{\kappa=1}^\mu N_\kappa'(K)$ (notice that $N'(K)$ is the number of cells intersected by more than one boundary component counted with multiplicity -- a cell intersected by two components is counted twice, a cell intersected by three components is counted three times etc.).

\subsection{The boundary decomposition}

We will now show that the Assumptions (R\ref{item:R_R2}) resp.\ (R\ref{item:R_Rd}) imply that the number of cells intersecting more than one boundary component $M_\kappa$ is small in a certain way. 

\begin{lemma}\label{l:N'bounded}
Let $K\subseteq \R{2}$ be a compact set satisfying (R\ref{item:R_R2}). 

Then $N'(rK+v) \le S$ for $v\in\mathbb{R}^2$ and sufficiently large $r$ for a bound $S$ which may depend on $K$, but not on $r$ or $v$. 
\end{lemma}
\prf Fix $\kappa\in\{1, \dots, \mu\}$. Since the functions $f_k, k=1,\dots, m,$ of (R\ref{item:R_R2}) are assumed to be either convex or concave, the set $M_\kappa$ is connected. Cells of $\mathcal{C}_\kappa'(rK+v)$ intersect also another boundary part $rM_{\kappa'}+v$, $\kappa'\ne\kappa$. By the compactness of the sets $M_\kappa$ there is $r_0$ such that for $r>r_0$ in such a situation always $rM_{\kappa}+v$ and $rM_{\kappa'}+v$ intersect. They intersect usually in one point, in some exceptional cases in two points. We will assume that there is only one intersection point in the following, since in the case of two intersection points the notation is blown up, while the ideas of the proof remain the same. 

If the angle $\alpha_{\kappa,\kappa'}$ between $M_{\kappa}$ and $M_{\kappa'}$ at their intersection point $z$ is bigger than $\pi/2$, then $p(\kappa)\ne p(\kappa')$ and the angle between two vectors from $G_{p(\kappa)}$ and from $G_{p(\kappa')}$ is at least $\pi/4$. Hence any cell intersecting both $rM_{\kappa}+v$ and $rM_{\kappa'}+v$ must contain a point which has distance at most $2(n-1)$ from $rz$ and therefore there can be at most $25n^2$ such cells (it would not be difficult to obtain a far lower bound; however, it only matters that this bound is independent of $r$, so we will not take the effort of improving it). 
 
So assume $\alpha_{\kappa,\kappa'}\le \pi/2$ from now on. Let $u^{(1)}$ and $u^{(2)}$ be the unit vectors such that $\pm u^{(1)}$ are normal vectors of $M_{\kappa}$ in $z$ and $\pm u^{(2)}$ are normal vectors of $M_{\kappa'}$ in $z$, oriented in such a way that both $u^{(1)}$ and $u^{(2)}$ point from $M_{\kappa'}$ to $M_{\kappa}$ in a neighborhood of $z$ (by the assumptions made so far, the angle which $M_{\kappa}$ and $M_{\kappa'}$ form at $z$ is strictly positive but not larger than $\pi/2$, so this choice is properly defined; it is convenient to orient the vectors like this and ignore whether they are now outward or inward normal vectors).  Let $\tilde u^{(1)}= \tfrac{2}{3} u^{(1)}+ \tfrac{1}{3} u^{(2)}$ and $\tilde u^{(2)}= \tfrac{1}{3} u^{(1)}+ \tfrac{2}{3} u^{(2)}$ and put $H_1:=\{y \in \R{2} \mid \langle y, \tilde u^{(1)} \rangle \ge\langle z, \tilde u^{(1)} \rangle \}$ and $H_2:=\{y \in \R{2} \mid \langle y, \tilde u^{(2)} \rangle \le\langle z, \tilde u^{(2)} \rangle \}$.  
There is some $\delta>0$ with $M_{\kappa} \cap B_\delta(z) \subseteq H_1$ and $M_{\kappa'} \cap B_\delta(z)\subseteq H_2$, where $B_\delta(z):=\{y\in\mathbb{R}^d \mid \|y-z\| < \delta \}$. For sufficiently large $r$ the sets $r(M_{\kappa} \setminus B_\delta(z))$ and $r(M_{\kappa'} \setminus B_\delta(z))$ have distance more than $(n-1)\sqrt{2}$ and hence there can be no cell which intersects both sets.
Let $\tilde\alpha_{\kappa,\kappa'}$ denote the angle between $H_1$ and $H_2$ and let $H$ denote a half-plane with $z\in \bd H$ and $M_{\kappa}\cup M_{\kappa'}\subseteq H$. Then a point of $H_1\cap H$ which lies in the same cell as a point from $H_2$  can have at most distance $(n-1)\sqrt{2}/\sin(\tilde\alpha_{\kappa,\kappa'})$ from $rz$. So there can be at most $((n-1)2\sqrt{2}/\sin(\tilde\alpha_{\kappa,\kappa'})+n)^2$ cells which intersect both $M_{\kappa}$ and $M_{\kappa'}$.  

\begin{center}
\begin{tikzpicture}
	\draw[fill] (0,0) circle(2pt);
	\draw (-0.2, -0.2) node{$z$};
	\draw (0.1, 2) node{$M_{\kappa}$};
	\draw (1.7, -0.2) node{$M_{\kappa'}$};
	\draw (2.4, 0.9) node{$H_2$};
	\draw (0.9, 2.4) node{$H_1$};
	\draw (-1.2,1.5) node{$H$};
	\draw[thick] (0.536,2) arc (150:180:4);
	\draw[thick] (2, -0.536) arc (60:90:4);
	\draw (0,0)--(2.6,1.3);
	\draw (0,0)--(1.3,2.6);
	\draw (-1.5,1.5)--(1,-1);
\end{tikzpicture}
\end{center}

Altogether we have 
\[ N'_\kappa(rK+v) \le \sum_{\kappa': \alpha_{\kappa,\kappa'}>\pi/2} 25n^2 + \sum_{\kappa': \alpha_{\kappa,\kappa'} \le \pi/2} ((n-1)2\sqrt{2}/\sin(\tilde\alpha_{\kappa,\kappa'})+n)^2, \]
where the sums are taken over all $\kappa' \in \{1,\dots,\mu\}$ with $M_\kappa \cap M_{\kappa'} \ne \emptyset$ (if $M_\kappa \cap M_{\kappa'}$ consists of two points, then $\kappa'$ contributes two summands to the sum).  
Summing up we get  
\[ N'(rK+v) \le \sum_{\kappa,\kappa': \alpha_{\kappa,\kappa'}>\pi/2} 25n^2 + \sum_{\kappa, \kappa': \alpha_{\kappa,\kappa'} \le \pi/2} ((n-1)2\sqrt{2}/\sin(\tilde\alpha_{\kappa,\kappa'})+n)^2, \]
where the sums are taken over all ordered pairs $(\kappa, \kappa') \in \{1,\dots, \mu \}^2$ with $M_\kappa \cap M_{\kappa'} \ne \emptyset$. \qed

\begin{lemma}\label{T:N'bounded}
Let $K\subseteq \R{d}$ be a compact set fulfilling (R\ref{item:R_Rd}). 

Then $N'(rK+v) \le Sr^{d-2}$ for $v\in\mathbb{R}^d$ and large enough $r$ for a bound $S$ which may depend on $K$, but not on $r$ or $v$. 
\end{lemma}

In the proof of this lemma we need the following lemma. Let $\tilde \kappa_{\tilde n}$ denote the Lebesgue measure of the $\tilde n$-dimensional unit ball.

\begin{lemma}\label{l:rect_cell} Let $R\subseteq\mathbb{R}^d$ be a rectangle of side-lengths $a_1,a_2,\dots, a_{d-1}$ within a hyperplane $H$. Assume $a_1\ge a_2 \ge \dots \ge a_{d-1}$. Then $R$ intersects at most $\sum_{\mathfrak{j}=0}^{d-1} \tilde \kappa_{d-\mathfrak{j}} \binom{d-1}{\mathfrak{j}} (n-1)^{d-\mathfrak{j}} d^{(d-\mathfrak{j})/2} a_1\cdot \dots \cdot a_\mathfrak{j} $ cells. 
\end{lemma}
\prf Consider the parallel set
\[ R_{\oplus \rho}:= \{ x\in\mathbb{R}^d \mid \|x-y\| \le \rho \mbox{ for one } y\in R\} , \quad \rho \ge 0,\]
of  $R$. By the Steiner formula its Lebesgue measure is given by
\[ \lambda_d( R_{\oplus \rho}) = \sum_{\mathfrak{j}=0}^d \tilde\kappa_{d-\mathfrak{j}} \rho^{d-\mathfrak{j}} V_\mathfrak{j}(R), \]
 where $V_\mathfrak{j}(R)$ is the $\mathfrak{j}$-th intrinsic volume of $R$; see e.g.\ \cite[(4.1)]{Schn14}. The intrinsic volumes of the rectangle $R$ are given by
\[ V_\mathfrak{j}(R) = \sum_{1\le i_\mathfrak{j}\le \dots \le i_\mathfrak{j}\le d-1} a_{i_1} \cdot \dots \cdot a_{i_\mathfrak{j}} \le \binom{d-1}{\mathfrak{j}} a_1\cdot a_2 \cdot\dots \cdot a_\mathfrak{j}, \quad \mathfrak{j}=0,\dots, d-1 \quad \mbox{ and } V_d(R)=0. \] 
Hence
\[ \lambda_d(R_{\oplus (n-1)\sqrt{d}}) \le \sum_{\mathfrak{j}=0}^{d-1} \tilde \kappa_{d-\mathfrak{j}} \binom{d-1}{\mathfrak{j}} (n-1)^{d-\mathfrak{j}} d^{(d-\mathfrak{j})/2} a_1\dots a_\mathfrak{j}. \]
A cell intersected by $R$ is completely covered by $R_{\oplus (n-1)\sqrt{d}}$. In particular, for a cell $[l_1,l_1+n-1] \times \dots \times [l_d, l_d+n-1]$ intersected by $R$, the subset $[l_1,l_1+1) \times \dots \times [l_d, l_d+1)$ is completely covered by $R_{\oplus (n-1)\sqrt{d}}$, and these subsets are disjoint for different cells. Hence the assertion follows. \qed\medskip

\prf[ of Lemma \ref{T:N'bounded}]
Fix $\kappa\in\{1,\dots, \mu\}$. Each cell of $\mathcal{C}_\kappa'(rK+v)$ intersects another boundary part $rM_{\kappa'}+v$, $\kappa\ne\kappa'$. By the compactness of the sets $M_\kappa$ there is $r_0$ such that for $r>r_0$ in such a situation $rM_{\kappa}+v$ and $rM_{\kappa'}+v$ always intersect. We have to distinguish several cases:

1.\ case: $M_{\kappa}$ and $M_{\kappa'}$ are parts of polytopes:\\
We may assume w.l.o.g.\ that $M_\kappa$ and $M_{\kappa'}$ are contained in hyperplanes. Since the angle under which $M_{\kappa}$ and $M_{\kappa'}$ meet is non-zero, their intersection is at most $(d-2)$-dimensional. 

Let $\tilde M$ be the set of points in $rM_{\kappa}+v$ that lie in a cell which is also intersected by $rM_{\kappa'}+v$. A point in $\tilde M$ can at most have distance $(n-1)\sqrt{d}$ from $rM_{\kappa'}+v$ and therefore it can most have distance $(n-1)\sqrt{d}/\sin(\alpha)$ from the affine hull $E$ of $M_\kappa \cap M_{\kappa'}$, where $\alpha$ is the angle under which $M_{\kappa}$ and $M_{\kappa'}$ intersect. However, the metric projection $p(E,x)$ of a point $x\in M_\kappa$ onto $E$ does not need not to lie in  $E \cap (M_\kappa \cup M_\kappa')$ and so  we have to find an upper bound for the diameter of $\{ p(E,x) \mid x\in \tilde M\}$, where $p(E,x)$ denotes the metric projection of $x$ onto $E$. For $\lambda\in\{\kappa, \kappa'\}$ consider the boundary $I_\lambda^{ + \rho}$ of the parallel set of $I_\lambda:=E\cap M_\lambda$ at distance $\rho$ within $E$. Let $\rho$ be small enough that $p(E,\mathfrak{v})$  either lies in $I_\lambda$ or has distance at least $\rho$ from $I_\lambda$ for any vertex $\mathfrak{v}$ of $M_\lambda$. Then
\[ \beta_\lambda:=  \min\big\{  d(M_\lambda \cap (x+E^\perp), x) \mid x\in I_\lambda^{ + \rho} \big\}/\rho > 0 , \]
where $d(M,x):=\inf\{ \| x-y\| \mid y\in M\}$, does not depend on $\rho$. Put $\beta:=\min\{\beta_\kappa, \beta_{\kappa'}\}$. 
 The diameter of $\{ p(E,x) \mid x\in \tilde M\}$ is at most $rl+2(n-1)\sqrt{d}(1 + 1/\beta/\sin(\alpha))$, where $l$ is the diameter of $I_{\kappa} \cup I_{\kappa'}$; indeed, since a point $x\in \tilde M$ can have distance at most $(n-1)\sqrt{d}/\sin(\alpha)$ from $E$, the point $p(E,x)$ can have distance at most $(n-1)\sqrt{d}(1+1/\sin(\alpha)/ \beta)$ from $I_{\kappa} \cup I_{\kappa'}$. 

Altogether, $\tilde M$ is contained in a $(d-1)$-dimensional rectangle with $d-2$ side lengths at most $rl+2(n-1)\sqrt{d}/ \beta/\sin(\alpha)$ and the remaining side length being at most $(n-1)\sqrt{d}/\sin(\alpha)$. Thus, by Lemma \ref{l:rect_cell}, the number of cells that are intersected both by $rM_{\kappa}+v$ and by $rM_{\kappa'}+v$ is bounded by a polynomial of degree $d-2$.


2.\ case: $M_{\kappa}$ and $M_{\kappa'}$ belong to the same smooth body $L$:\\
The \emph{support function} of a non-empty compact set $\tilde L\subseteq \mathbb{R}^d$ is defined as 
\[ h_{\tilde L}(u) := \max\{ \langle x, u\rangle \mid x \in \tilde L\}, \quad u \in \mathbb{R}^d; \]
see \cite[Sec.\ 1.7.1]{Schn14}. From \cite[p.\ 115]{Schn14} we get that $h_L$ is twice differentiable on $\mathbb{R}^d\setminus \{0\}$, since $L$ is convex and $\bd L$ is a $C^2$-manifold with non-vanishing Gauss-Kronecker curvature. Let $H_uh_L$ be its Hessian matrix in $u\in\mathbb{R}^d\setminus\{0\}$ and put $Li:= \max\{ \|H_uh_L\|_2 \mid u\in S^{d-1} \}$, where $\| \cdot \|_2$ denotes the matrix norm induced by the Euclidean norm. 

The sets $M_{\kappa}$ and $M_{\kappa'}$ belong to two different sets $H_{p}$ and $H_{p'}$ for $p,p'\in S(n^d)$. In any point $z\in M_\kappa \cap M_{\kappa'}$ there must be a normal vector of $L$ that lies in $\tilde G:=G_{p} \cap G_{p'}\subseteq S^{d-1}$. By \cite[Corollary 1.7.3]{Schn14} $L$ can have an exterior normal vector in $\tilde G$ only in points that lie in the image $\tilde H$ of $\tilde G$ under $\nabla h_{rL+v}$. The set $\tilde G$ is a subset of a $(d-2)$-dimensional sphere. Hence it can be covered by $2(d-1)$ sets isometric to $\{(\tilde u_1, \dots, \tilde u_{d-1}) \in S^{d-2}\mid \tilde u_{d-1}\ge 1/\sqrt{d-1}\}$. This set is the image of $\tilde D:=\{ \tilde v\in\mathbb{R}^{d-2} \mid \|\tilde v\|\le \sqrt{(d-2)/(d-1)}\}$ under the mapping $\tilde f:(\tilde v_1,\dots, \tilde v_{d-2}) \mapsto (\tilde v_1,\dots, \tilde v_{d-2}, \sqrt{1-\tilde v_1^2-\dots-\tilde v_{d-2}^2}) $ which is Lipschitz continuous with Lipschitz constant $\sqrt{d}$. There are $\tilde N:=(\lceil d\cdot r\cdot Li\rceil)^{d-2}$ points $\tilde v^{(1)}, \dots, \tilde v^{(\tilde N)}$ with $\tilde D\subseteq \bigcup_{\mathfrak{k}=1}^{\tilde N}B_{1/(\sqrt{d}\cdot r\cdot Li)}(\tilde v^{(\mathfrak{k})})$ and thus $\tilde f(\tilde D) \subseteq_{\mathfrak{k}=1}^{\tilde N} B_{1/(r\cdot Li)}(\tilde f(\tilde v_\mathfrak{k}))$, where $\lceil a \rceil := \max\{b\in \mathbb{Z} \mid b\le a \}$. Hence there are $N:=2(d-1)\cdot \tilde N=2(d-1) \cdot (\lceil d\cdot r\cdot Li\rceil)^{d-2}$ points $w^{(1)}, \dots, w^{(N)} \in \tilde G$ such that $\tilde G\subseteq \bigcup_{\mathfrak{k}=1}^N B_{1/(r\cdot Li)}(w^{(\mathfrak{k})})$. Now for every point $z\in \tilde H$ one of the points $x^{(\mathfrak{k})}:=\nabla h_{rL+v}(w^{(\mathfrak{k})})$, $\mathfrak{k}=1, \dots, N$ has distance at most $1$, since $\nabla h_{rL+v}$ has Lipschitz constant $r\cdot Li$.

Since the principle curvatures depend continuously on the point, a compactness argument ensures that the principle radii of curvature of $rL$ are bounded from below by $(n-1)\sqrt{d}$ for sufficiently large $r$. Then a point in a cell intersecting both $rM_{\kappa}+v$ and $rM_{\kappa'}+v$ can have distance at most $(n-1)\sqrt{d}$ from the nearest point in $\tilde H$ and therefore it has distance less than $n\sqrt{d}$ from the nearest point $x^{(\mathfrak{k})}, \mathfrak{k}\in\{1,\dots, N\}$. Thus there can be at most $(n+2n\sqrt{d})^d \cdot N = 2(d-1)(n+2n\sqrt{d})^d \cdot (\lceil d\cdot r\cdot Li\rceil)^{d-2}$ cells intersecting both $rM_{\kappa}+v$ and $rM_{\kappa'}+v$. 

3.\ case: $M_{\kappa}$ belongs to a polytope $P$, while $M_{\kappa'}$ belongs to a smooth convex body $L$ (or the other way round):\\
Let $E_1, \dots, E_s$ denote the affine hulls of the facettes of $P$ which are intersected by $M_{\kappa} \cap M_{\kappa'}$. Then $\Gamma^{(i)}:=E_i \cap {\bd L}$ is the boundary of a convex body lying in $E_i$ for each $i=1,\dots, s$. Put $\tilde \Gamma^{(i)} := \Gamma^{(i)} \cap M_{\kappa} \cap M_{\kappa'}$. 

By the smoothness assumption on $L$, the angle $E_i$ and $M_{\kappa'}$ form at the points $z \in \Gamma^{(i)}$ is continuous as function of $z$. By the compactness of $\tilde \Gamma^{(i)}$ it attains a minimum $\alpha_i>0$. Unlike in the 1st case one cannot assume that two points $x\in M_\kappa$ and $x'\in M_{\kappa'}$ are always seen from an appropriate point $z\in\Gamma^{(i)}$ under an angle of at least $\alpha$, since $\bd L$ is curved. We shall explain now, why this can be assumed with $\alpha$ replaced by $\alpha/2$. Since $\Gamma^{(i)}$ is a $C^2$-manifold, there is  some critical radius $\rho_*$ such that the metric projection $p(\Gamma^{(i)},x)$ is defined uniquely for any point $x$ of distance less than $\rho_*$ to $\Gamma^{(i)}$. Let $e(z)$ for $z\in \Gamma^{(i)}$ be the infimum over all $\rho>0$ for which there is $\tau\in (0,\rho\cdot \sin{(\alpha/2)})$ such that one of the four points $z \pm\rho \nu(z) \pm\tau \nu_{E_i}$ lies in $\bd L$, where $\nu(z)$ is a unit normal vector of $\Gamma^{(i)}$ in $z\in \Gamma^{(i)}$ within the linear subspace which is parallel to $E_i$ and where $\nu_{E_i}$ is the unit normal vector of $E_i$. Now $e$ is lower semicontinuous, since if $(z_n)_{n\in\mathbb{N}}$, $(\rho_n)_{n\in\mathbb{N}}$ and $(\tau_n)_{n\in\mathbb{N}}$ with $z_n+\rho_n \nu(z_n)+ \tau_n \nu_{E_i} \in \bd L$ for all $n\in\mathbb{N}$ converge to limits $z$, $\rho$ and $\tau$, then $z+\rho \nu(z)+ \tau \nu_{E_i} \in \bd L$. Hence $e$ attains a minimum $e_*>0$ on $\tilde \Gamma^{(i)}$. Let $\epsilon>0$ be such that $e(z)>e_*/2$ for any $z\in \Gamma^{(i)} \cap \tilde \Gamma^{(i)}_{\oplus \epsilon}$. Let $r$ be large enough such that the distance of $rE_i\setminus r\tilde\Gamma^{(i)}_{\oplus \min\{\rho_*, e_*, \epsilon\}/2}$ to $r\bd L$ is larger than $(n-1)\sqrt{d}$. Then cell intersecting both $(rE_i+v)\cap (rM_{\kappa}+v)$ and $rM_{\kappa'}+v$ must contain a point of distance at most $(n-1)\sqrt{d}/\sin(\alpha/2)$ from $\Gamma^{(i)}$. 

Now $\Gamma^{(i)}$ can be represented as union of the graphs of $2(d-1)$ convex function, which are Lipschitz continuous with Lipschitz constant $1$. Hence there are $N^{(i)}:= (2d-2)\cdot (\lceil \sqrt{(d-2)/2}\cdot r\cdot \Lambda^{(i)}\rceil)^{d-2}$ points $x_1^{(i)},\dots, x_{N^{(i)}} ^{(i)}\in r\Gamma^{(i)}$ with $\bigcup_{\mathfrak{k}=1}^{N^{(i)}} B_1(x_\mathfrak{k}^{(i)}) \supseteq r\Gamma^{(i)}$, where $\Lambda^{(i)}$ is the diameter of $\Gamma^{(i)}$. Now a cell intersecting both $(rE_i+v)\cap (rM_{\kappa}+v)$ and $rM_{\kappa'}+v$ must contain a point of distance at most $(n-1)\sqrt{d}/\sin(\alpha/2)+1$ from the nearest point $x_\mathfrak{k}^{(i)},\, \mathfrak{k}=1,\dots, N^{(i)}.$ Hence there are less than 
\[ \sum_{i=1}^s N^{(i)}\cdot (\lceil 2(n-1)\sqrt{d}/\sin(\alpha/2)+n+2\rceil)^d = \sum_{i=1}^s (2d-2)\cdot (\lceil \sqrt{(d-2)/2}\cdot r\cdot \Lambda^{(i)}\rceil)^{d-2}\cdot (\lceil 2(n-1)\sqrt{d}/\sin(\alpha/2)+n+2\rceil)^d\] 
cells intersecting both $M_\kappa$ and $M_{\kappa'}$.

4.\ case: $M_{\kappa}$ and $M_{\kappa'}$ belong to different smooth bodies $L_{k_1}$ and $L_{k_2}$:\\ 
Fix $z_0\in \bd L_{k_1} \cap \bd L_{k_2}$. Then there are a neighborhood $U\subseteq \R{d}$ of $z_0$, a vector $\nu_1\in S^{d-1}$, a neighborhood $W$ of $0$ in $\nu_1^\perp$ and two $C^2$-functions $g_1, g_2: W\to \mathbb{R}$ such that 
\[ \bd L_{k_i} \cap U = \{w + g_i(w) \nu_1 + z_0 \mid w \in W \}, \quad i=1,2. \]
By the implicit function theorem there is a unit vector $\nu_2 \in W$, a neighborhood $W'$ of $0$ in $W\cap \nu_2^\perp$ and a $C^2$-function $h: W'\to \mathbb{R}$ with  
\[ \Gamma := \bd L_{k_1}\cap \bd L_{k_2} \cap U =\{ w + h(w)\nu_2 + g_1\big(w+h(w)\nu_2\big)\nu_1 +z_0\mid w\in W'\} \]
possibly after replacing $U$ by a smaller set.

 Replacing $W'$ by a subset if necessary we may assume that $g_1$, $g_2$ and $h$ are Lipschitz continuous with Lipschitz constant $1$. Moreover, choose $\Lambda$ such that $W'$ is contained in a cube of side-length $\Lambda$. Then there are $N:= (\lceil\Lambda\sqrt{d-2}\cdot r\rceil)^{d-2}$ points $x_1,\dots, x_N \in \Gamma$ with $\bigcup_{\mathfrak{k}=1}^N B_{1/r}(x_\mathfrak{k}) \supseteq \Gamma$, namely the images of the points of the lattice covering $W'$ with grid distance $\tfrac{1}{\sqrt{d-2}r}$ under the mapping $w\mapsto w + h(w)\nu_2 + g_1\big(w+h(w)\nu_2\big)\nu_1  +z_0$. By the compactness of $\bd L_{k_1}\cap \bd L_{k_2}$ there are finitely many of the manifolds $\Gamma^{(1)}, \dots, \Gamma^{(s)}$ constructed above such that $\bd L_{k_1}\cap \bd L_{k_2} = \bigcup_{i=1}^s \Gamma^{(i)}$. Let $N^{(i)}$, $Li^{(i)}$, $\Lambda^{(i)}$ denote the numbers constructed above associated to the manifold $\Gamma^{(i)}$, $i=1,\dots, s$ -- notice that the $N^{(i)}$ depend on $r$, while the $Li^{(i)}$ and $\Lambda^{(i)}$ are independent of $r$. So there are $N^\Sigma:=\sum_{i=1}^s N^{(i)}$ points $x_1,\dots, x_{N^\Sigma}\in r\bd L_{k_1}\cap r\bd L_{k_2}$ with $\bigcup_{\mathfrak{k}=1}^{N^\Sigma} B_1(x_\mathfrak{k}) \supseteq r\bd L_{k_1}\cap r\bd L_{k_2}$. 

Similar as in the 3rd case we get that there is a minimal angle $\alpha$ under which $\bd L_{k_1}$ and $\bd L_{k_2}$ intersect. Now fix some $i\in\{1, \dots, s \}$ and let $e(z),\, z\in \Gamma^{(i)}$, be the infimum over all $\rho>0$ for which there are $\tau\in\mathbb{R}$ such that one of the four points $z\pm \rho \nu(z)\pm \tau \nu_1$ lies in $\bd L_{k_1}$ and satisfies $d(\bd L_{k_2},z\pm \rho \nu(z)\pm \tau \nu_1)\le \sqrt{\rho^2+\tau^2}\cdot \sin(\alpha/2)$, where $\nu(z)$ is the normal vector of $\{w+h(w)\nu_2 \mid w\in W'\}$ at $z$ within $\nu_1^\perp$. Similar as in the 3rd case it can be shown at $e$ attains its minimum.

The same way as in the 3rd case one sees that for sufficiently large $r$ there are at most 
\[\sum_{i=1}^s N^{(i)}(\lceil 2(n-1)\sqrt{d}/\sin(\alpha/2)+n+2\rceil)^d = \sum_{i=1}^s (\lceil\Lambda^{(i)}\sqrt{d-2}\cdot r\rceil)^{d-2}(\lceil 2(n-1)\sqrt{d}/\sin(\alpha/2)+n+2\rceil)^d \] 
cells intersecting both $rM_{\kappa}+v$ and $rM_{\kappa'}+v$.\medskip

So for any pair $(\kappa, \kappa')\in \{1, \dots, \mu\}$ with $\kappa \ne \kappa'$ the number of cells intersecting both $rM_{\kappa}+v$ and $rM_{\kappa'}+v$ is bounded uniformly in $v$ by a function of order $d-2$ in $r$. Summing up over all such pairs we get that $N'(rK+v)$ is bounded uniformly in $v$ by a function of order $d-2$ in $r$. \qed

\subsection{Bounds on the pixel configuration counts}

Our aim in this section is to derive bounds on the numbers $N_{\kappa,j}^-$ defined in Section \ref{ss:notation}. 

Let $K\subseteq \mathbb{R}^d$ be a compact set fulfilling (R\ref{item:R_R2}) if $d=2$ or (R\ref{item:R_Rd}) if $d>2$. For every $p\in S(n^d)$ and $j=1,\dots, 2^{(n^d)}$, there are $\epsilon^-=(\epsilon_{1,j,p}^-, \dots, \epsilon_{d,j,p}^-), \epsilon^+ =(\epsilon_{1,j,p}^+, \dots, \epsilon_{d,j,p}^+) \in \{ 0,\dots, n-1\}^d$ with
\[ l+ B_j \subseteq K,\ l+W_j\subseteq K^C  \iff  l+\epsilon^- \in K,\ l+\epsilon^+ \in K^C, \]
provided that $\bd K \cap (l+[0,(n-1)]^d) \subseteq H_p$, $l\in\mathbb{Z}^d$. We get
\begin{align*}
 N_{\kappa,j}^- &=\#\{l\in\mathbb{Z}^d \mid  l+\epsilon^- \in K, \, l+\epsilon^+ \in K^C, \emptyset \ne (\bd K) \cap (l+[0,n-1]^d) \subseteq M_\kappa \}\\
& = \sum_{i=1}^d \sum_{\gamma=\epsilon^-_{i,j,p(\kappa)}}^{\epsilon^+_{i,j,p(\kappa)}-1} 
\begin{aligned}[t]\#\{l\in\mathbb{Z}^d \mid& l+ (\epsilon_{1,j,p(\kappa)}^-, \dots, \epsilon_{i-1,j,p(\kappa)}^-, \gamma, \epsilon_{i+1, j,p(\kappa)}^+, \dots, \epsilon_{d, j,p(\kappa)}^+) \in K, \\
& l+ (\epsilon_{1,j,p(\kappa)}^-, \dots, \epsilon_{i-1,j,p(\kappa)}^-, \gamma+1, \epsilon_{i+1, j,p(\kappa)}^+, \dots, \epsilon_{d, j,p(\kappa)}^+) \in K^C,\\
& \emptyset \ne (\bd K) \cap (l+[0,n-1]^d) \subseteq M_\kappa \}\\
-\#\{l\in\mathbb{Z}^d \mid& l+ (\epsilon_{1,j,p(\kappa)}^-, \dots, \epsilon_{i-1,j,p(\kappa)}^-, \gamma+1, \epsilon_{i+1, j,p(\kappa)}^+, \dots, \epsilon_{d, j,p(\kappa)}^+) \in K, \\
& l+ (\epsilon_{1,j,p(\kappa)}^-, \dots, \epsilon_{i-1,j,p(\kappa)}^-, \gamma, \epsilon_{i+1, j,p(\kappa)}^+, \dots, \epsilon_{d, j,p(\kappa)}^+) \in K^C,\\
& \emptyset \ne (\bd K) \cap (l+[0,n-1]^d) \subseteq M_\kappa \}
\end{aligned}
\end{align*}
(in the second set $\gamma$ and $\gamma+1$ are exchanged), where we understand $\sum_{\gamma=b}^a \alpha_\gamma= -\sum_{\gamma=a+1}^{b-1} \alpha_\gamma$ if $a+1<b$ and $\sum_{\gamma=a}^{a-1} \alpha_\gamma=0$. 

Thus we introduce the following notation: 
If $x=(x_1, \dots, x_{d-1})\in\R{d-1}$, $y\in\R{}$ and $i\in\{1,\dots, d\}$, then $[x]^i_y:=(x_1, \dots, x_{i-1}, y, x_{i}, \dots, x_{d-1})$ is the vector in $\R{d}$ obtained from $x$ by inserting $y$ at the $i$-th position and moving the $i$-th to $(d-1)$-st entry by one position. 

For $\kappa\in\{1,\dots, \mu\}$, $\epsilon\in \{0,\dots,n-1\}^{d-1}$, $\gamma\in\{0,\dots, n-2\}$ and $i\in\{1,\dots, d\}$ let $\mathcal{C}_{\kappa,\epsilon,\gamma,i}$ denote the system of cells which intersect $\bd K$, but only in points belonging to $M_\kappa$, and for which $\bd K$ passes through between $l+[\epsilon]^i_\gamma$ and $l+ [\epsilon]^i_{\gamma+1}$. Formally, put
\[
 \mathcal{C}_{\kappa,\epsilon,\gamma,i} :=
 \{l\in \mathbb{Z}^d \mid  (l+[0,n-1]^d) \cap \bd K \neq \emptyset, \ (l+[0,n-1]^d) \cap \bd K \subseteq M_\kappa, \
  \# (K \cap \{l+[\epsilon]^i_\gamma, l+[\epsilon]^i_{\gamma+1} \} ) = 1 \} , 
\]  
and $N^-_{\kappa,\epsilon,\gamma,i} := \# \mathcal{C}_{\kappa,\epsilon,\gamma,i}$; see Figure \ref{F:C}. 

Now we have
\begin{equation}
N_{\kappa,j}^- = \sum_{i=1}^d \sum_{\gamma=\epsilon^-_{i,j,p(\kappa)}}^{\epsilon^+_{i,j,p(\kappa)}-1} \sign(u_i)\cdot N^-_{\kappa, (\epsilon_{1,j,p(\kappa)}^-, \dots, \epsilon_{i-1,j,p(\kappa)}^-, \epsilon_{i+1, j,p(\kappa)}^+, \dots, \epsilon_{d, j,p(\kappa)}^+ ), \gamma,i}
\label{e:N=N}\end{equation}
for an arbitrary $(u_1, \dots, u_d)\in\tilde G_{p(\kappa)}$, where $\sign$ is the sign function. 


\begin{figure}%
\begin{center}
\begin{tikzpicture}
\foreach \x in {0,0.5,1, 2.5,3,3.5, 5,5.5,6}
	\foreach \y in {0,0.5,1}
	\draw (\x,\y) circle(2pt);
\draw (-0.5,1.2)--(1.5,0);
\draw (	2, 1)--(4,-0.2);
\draw (4.5,0.84)--(6.5,-0.34);
\foreach \v in {(0.0,0.5), (0.5,0.5), (0.0,0.0), (0.5,0.0), (1.0,0.0), (2.5,0.5), (2.5,0.0), (3.0,0.0), (3.5,0.0), (5.0,0.0), (5.0,0.5), (5.5,0.0)}
\draw[fill] \v circle(2pt); 
\end{tikzpicture}
\end{center}
The cells having one of these three pixel configurations make up the set $\mathcal{C}_{\kappa,\epsilon,\gamma,i}$ in the case $d=2$, $n=3$, $\epsilon=0$, $\gamma=1$, $i=2$ and $\kappa\in \{1,\dots,\mu\}$ is as indicated by the slope of the line (as usual, the first coordinate of a point is plotted horizontal, increasing to the right, and the second coordinate is plotted vertical, increasing upwards). The defining property is that the border line must pass through between the upper left point and the middle left point.
\caption{An example of a set $\mathcal{C}_{\kappa,\epsilon,\gamma,i}$}%
\label{F:C}%
\end{figure}

We are going to show that for $\kappa\in\{1, \dots, \mu\}$, $\epsilon\in \{0,\dots,n-1\}^{d-1}$, $\gamma\in\{0,\dots, n-2\}$ and $i\in\{1,\dots,d\}$ the number $N^-_{\kappa,\epsilon,\gamma,i}$ approximately equals
\[ I_{\kappa,i} := \lambda_{d-1}(M_\kappa\mid E_i),\]
where $M|E$ denotes the image of $M$ under the orthogonal projection onto $E$ and $E_i:=e_i^\perp$ for the standard base $e_1, \dots, e_d$ of $\mathbb{R}^d$.

Before we come to the main result of this subsection, we discuss the special case that $d=2$, $n=2$ and $p\in S(4)$ is such that $G_p=\{ (u_1, u_2) \in\mathbb{R}^2 \mid 0 \le u_1 \le u_2 \}$. 

In this case $I_{\kappa,1}$ is the length of the orthogonal projection of $M_\kappa$ onto the $y$-axis and $I_{\kappa,2}$ is the length of the projection of $M_\kappa$ onto the $x$-axis. 
 
 We let in this special case $N_{\kappa,1}^-$ denote the number of cells $C$ with $\emptyset \ne C \cap(\bd K) \subseteq M_\kappa$ with a black pixel in the lower left corner and the remaining three pixels being white, $N_{\kappa,2}^-$ is the number of cells ``on'' $M_\kappa$  with two black pixels in the lower row and two white pixels in the upper row and $N_{\kappa,3}^-$ is the same number for the configuration in which only the pixel in the upper right corner is white. Thus we have
\begin{align*}
N_{\kappa,0,0,1}^- &= N_{\kappa,1}^-  								&  N_{\kappa,1,0,1}^- &= N_{\kappa,3}^-   \\   
N_{\kappa,0,0,2}^- &= N_{\kappa,1}^- + N_{\kappa,2}^- &  N_{\kappa,1,0,2}^- &=N_{\kappa,2}^- + N_{\kappa,3}^-.
\end{align*} 

\begin{theorem}
Let $K\subseteq \mathbb{R}^2$ fulfill (R\ref{item:R_R2}) and let $M_\kappa\subseteq H_p$ for the $p\in S(4)$ described above. Then
\begin{align}
 N^-_{\kappa,3} - N'_\kappa \le I_{\kappa,1}&					 		\mbox{ and }  I_{\kappa,1} \le N^-_{\kappa,3}+N'_\kappa, \label{e:bound3}\\
 N^-_{\kappa,1} - N'_\kappa \le I_{\kappa,1}&							\mbox{ and }  I_{\kappa,1} \le N^-_{\kappa,1}+N'_\kappa, \label{e:bound1}\\
N^-_{\kappa,3} + N^-_{\kappa,2} - N'_\kappa \le I_{\kappa,2} &	\mbox{ and }  I_{\kappa,2}  \le N^-_{\kappa,3} + N^-_{\kappa,2} + N'_\kappa, \label{e:bound23}\\
N^-_{\kappa,1} + N^-_{\kappa,2} - N'_\kappa \le I_{\kappa,2} &	\mbox{ and }  I_{\kappa,2}  \le N^-_{\kappa,1} + N^-_{\kappa,2} + N'_\kappa. \label{e:bound12}
\end{align}
\end{theorem}

\begin{figure}%

\begin{center}
\includegraphics[width=0.4\columnwidth]{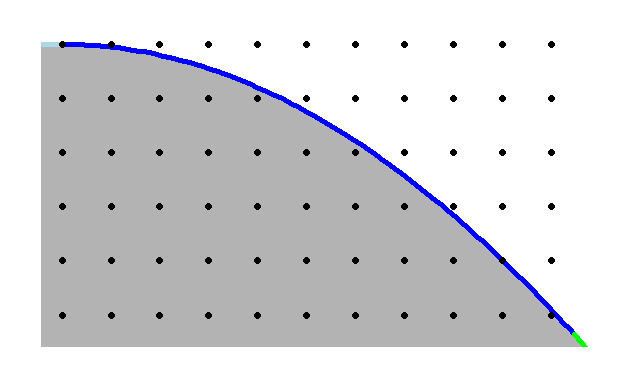}%
\end{center}
\caption{A component $M_\kappa$ of the boundary of $K$.}%
\label{F:component}%
\end{figure}

Since this theorem is a special case of Theorem \ref{T:pixel_bounds} below, we do not give a formal prove, but we only describe the ideas. A visualization is given in in Figure \ref{F:component}. 

A horizontal line intersecting $M_\kappa$ will pass through exactly one cell with three black pixels and one white pixel (unless the horizontal line is very close to one the endpoints of $M_\kappa$). Thus the number $N_{\kappa,3}^-$ of such cells equals the length of the projection of $M_\kappa$ onto the $y$-axis up to some ``border'' effects coming from the fact that we had to exclude horizontal lines which are close to one of the endpoints of $M_\kappa$ in the sentence before. If we make this approximation precise by two inequalities we arrive at \eqref{e:bound3}. 

Exactly the same argument is true when cells with three black pixels are replaced by cells with one black pixel. Hence \eqref{e:bound1} follows.

We see that there are vertical lines intersecting no cell with three black pixels and one white pixel even if the line intersects $M_\kappa$ far from the endpoints. However, we can observe that every vertical line either intersects a cell with three black pixels or a cell with two black pixels and that it can only intersect one such cell. Hence the sum of the numbers of these two types of cells equals the length of the orthogonal projection of $M_\kappa$ onto the $x$-axis. Similar as above we obtain \eqref{e:bound23}. The same way as \eqref{e:bound23} we get \eqref{e:bound12}.

Now let us return to arbitrary $n$ and $d$. 

\begin{theorem}\label{T:pixel_bounds}
Let $K\subseteq \R{d}$ fulfill (R\ref{item:R_R2}) resp.\ (R\ref{item:R_Rd}). Let $\kappa\in \{1,\dots,\mu\}$, $i\in\{1,\dots,d\}$, $\epsilon\in\{0,\dots, n-1\}^{d-1}$ and $\gamma\in\{0,\dots,n-2\}$. Then
\begin{enumerate}[(i)]
\item $N^-_{\kappa,\epsilon, \gamma,i} - N'_\kappa \le I_{\kappa,i}$  \label{item:Iungl}
\item $ I_{\kappa,i} \le N^-_{\kappa,\epsilon,\gamma,i} + N'_\kappa$ \label{item:Jungl}
\end{enumerate}
\end{theorem}
\prf 
To prove (\ref{item:Iungl}) let $\mathcal{Z} \subseteq \mathbb{Z}^{d-1}$ denote the set of points $l=(l_1,\dots, l_{d-1})\in\mathbb{Z}^{d-1}$ for which there is some $k\in\mathbb{Z}$ with 
\[ [l_1,l_1+n-1] \times \dots \times [l_{i-1},l_{i-1}+n-1]\times [k, k+n-1] \times [l_{i},l_{i}+n-1] \times \dots \times [l_{d-1},l_{d-1}+n-1] \in \mathcal{C}_{\kappa,\epsilon,\gamma,i}. \]
Clearly this integer $k$ is determined uniquely - denote it by $k(l)$ and put $C_l:=[l_1,l_1+1] \times \dots \times [l_{i-1},l_{i-1}+1]\times [k(l), k(l)+1] \times [l_{i},l_{i}+1] \times \dots \times [l_{d-1},l_{d-1}+1]$. Let $\mathcal{Z}_1:= \{l\in\mathcal{Z} \mid \times_{s} [l_s,l_s+1] \subseteq (M_\kappa \cap C_l)|E_i \}$, where $E_i:=e_i^\perp $, and denote $\mathcal{Z}_2:=\mathcal{Z} \setminus\mathcal{Z}_1$. Then 
\[ \# \mathcal{Z}_1 \le \lambda_{d-1}\big( \bigcup_{l\in\mathcal{Z}_1} (M_\kappa \cap C_l)| E_i \big) \le I_{\kappa,i}. \] 
For every $l\in\mathcal{Z}_2$ there is some integer $k'(l)$ with 
\[ [l_1,l_1+n-1] \times \dots \times [l_{i-1},l_{i-1}+n-1]\times [k'(l), k'(l)+n-1] \times [l_{i},l_{i}+n-1] \times \dots \times [l_{d-1},l_{d-1}+n-1] \in \mathcal{C}'_\kappa. \]
Thus $\# \mathcal{Z}_2 \le N'_\kappa$. Hence the first inequality follows.

In order to show (\ref{item:Jungl}), we identify $E_i$ with $\mathbb{R}^{d-1}$.
 Let $(l_1,\dots, l_{d-1})\in \mathbb{Z}^{d-1}$ be a point with $M_\kappa|E_i$ intersecting $\times_{s} [l_s,l_s+n-1]$. If $\times_{s} [l_s,l_s+n-1] \subseteq M_\kappa|E_{i}$ does not hold, then $(M_\kappa \cap M_{\kappa'}) |E_{i}$ intersects $\times_{s} [l_s,l_s+n-1]$ for a $\kappa'\ne \kappa$. Thus there is a cell with coordinates in $\times_{s} [l_s,l_s+n-1]$ that belongs to $\mathcal{C}'$. So assume $\times_{s} [l_s,l_s+n-1] \subseteq M_\kappa|E_{i}$ from now on. 

Let $k\in\mathbb{Z}$ be the integer for which the line segment from $[l+\epsilon]_{k+\gamma}^{i}$ to $[l+\epsilon]_{k+\gamma+1}^{i}$ intersects $M_\kappa$ - in case this $k$ is not unique, choose the largest one if the $i$-th component of vectors from $\tilde G_{p(\kappa)}$ is positive and choose the smallest one if the $i$-th component of the vectors from $\tilde G_{p(\kappa)}$ is negative.  
Obviously, 
\[ [l_1, l_1+n-1] \times \dots \times [k, k+n-1] \times \dots \times [l_{d-1}, l_{d-1}+n-1] \in\mathcal{C}_{\kappa,\epsilon,\gamma, i}\cup \mathcal{C}'_\kappa.\]

Thus there is either a cell contributing to $N^-_{\kappa,\epsilon,\gamma,i}$ or a cell contributing to $N'_\kappa$ with coordinates in $\times_{s\ne i} [l_s,l_s+n-1]$ unless $(M_\kappa|E_i) \cap (\times_{s\ne i} [l_s,l_s+n-1]) = \emptyset$. In the latter case, we have, however, $\lambda_{d-1}((M_\kappa|E_i) \cap (\times_{s\ne i} [l_s,l_s+1]) )=0$. Summing up over all $l\in\mathbb{Z}^{d-1}$ yields the second inequality. 
\qed \medskip

For $u\in S^{d-1}$ and $c\in \mathbb{R}$ we let $(B_{u,c}, W_{u,c})$ denote the pixel configuration with $B_{u,c}=\{x\in \{0,\dots,n-1\}^d \mid \langle x, u \rangle \le c \}$ and $W_{u,c}=\{x\in \{0,\dots,n-1\}^d \mid \langle x, u \rangle > c \}$. 

\begin{kor}\label{K:pixel_bounds}
Let $K\subseteq \R{d}$ fulfill (R\ref{item:R_R2}) resp.\ (R\ref{item:R_Rd}). Let $\kappa\in \{1,\dots,\mu\}$ and $j\in\{1,\dots,2^{n^d}\}$. If $(B_j,W_j)=(B_{u,c}, W_{u,c})$ for some $u\in \tilde G_{p(\kappa)}$ and $c\in\mathbb{R}$, then
\begin{gather*} 
\sum_{i=1}^d \big((\epsilon^+_{i,j,p(\kappa)} - \epsilon^-_{i,j,p(\kappa)})\cdot \sign(u_i) \cdot I_{\kappa,i} - |\epsilon^+_{i,j,p(\kappa)} - \epsilon^-_{i,j,p(\kappa)}|\cdot  N'_\kappa \big) \ \le \ N^-_{\kappa,j}\\
 \le \sum_{i=1}^d \big((\epsilon^+_{i,j,p(\kappa)} - \epsilon^-_{i,j,p(\kappa)})\cdot\sign(u_i) \cdot  I_{\kappa,i} + |\epsilon^+_{i,j,p(\kappa)} - \epsilon^-_{i,j,p(\kappa)}|\cdot  N'_\kappa \big).
\end{gather*}
If $(B_j,W_j)$ is not of this form, then
\[ N^-_{\kappa,j}=0.\]
\end{kor}
\prf The first assertion is an immediate consequence of \eqref{e:N=N} and Theorem \ref{T:pixel_bounds}, while the second one is trivial. \qed

\subsection{Proof of the main result}

Here we proof Theorem \ref{T:main}. 

The idea of this proof is the following: Using Corollary \ref{K:pixel_bounds} we can approximate $N_{\kappa,j}^-(\frac{1}{t}K+v)$ up to a small correction term by an expression  which is invariant under translations. 
For the translation invariant expression the supremum and the infimum on the left hand side of (\ref{e:main}) cancel. Hence we are left with the correction term, which grows at most of order $(1/t)^{d-2}$ by Lemma \ref{l:N'bounded} and Lemma \ref{T:N'bounded}. Together with the normalizing factor $t^{d-1}$ in \eqref{e:form_local} the whole expression thus tends linearly to zero as $t\to 0$.

\prf[ of Theorem \ref{T:main}] 
We have
\begin{align*}
\sum_{\kappa=1}^\mu N_{\kappa,j}^-(\tfrac{1}{t}K+v) \le  N_{t,j}(K+tv) & \le \sum_{\kappa=1}^\mu N_{\kappa,j}^-(\tfrac{1}{t}K+v) + N'(\tfrac{1}{t}K+v),\quad j=1,\dots, 2^{(n^d)}.
\end{align*} 


From Corollary \ref{K:pixel_bounds} we get
\begin{align*}
\sup\{& \hat S_{t}(K+tv) \mid v\in[0,1)^d \}  - \inf\{ \hat S_{t}(K+tv) \mid v\in[0,1)^d \} \\
&\le 	\begin{aligned}[t] \sum_{j=1}^{2^{(n^d)}}  		\big(   \sup\{ w_j\cdot t^{d-1} \cdot N_{t,j}(K+tv) \mid v\in[0,1)^d \} - \inf\{ w_j\cdot t^{d-1} \cdot N_{t,j}(K+tv) \mid v\in[0,1)^d \}		\big) \end{aligned}
\\
&\le 	\begin{aligned}[t] \sum_{j=1}^{2^{(n^d)}} \sum_{\kappa=1}^\mu & |w_j|\cdot t^{d-1} \cdot 		\big( \sup\{ N_{\kappa,j}^-(\tfrac{1}{t}K+v) \mid v\in[0,1)^d \} - \inf\{ N_{\kappa,j}^-(\tfrac{1}{t}K+v) \mid v\in[0,1)^d \}		\big) \\
									&+ \Big( \sum_{j=1}^{2^{(n^d)}} |w_j|\Big) \cdot t^{d-1} \cdot  \sup\{N'(\tfrac{1}{t}K+v) \mid v\in[0,1)^d \}\end{aligned}
\\
&\le \begin{aligned}[t] \sum_{j=1}^{2^{(n^d)}} &\sum_{\kappa=1}^\mu |w_j|\cdot t^{d-1} \cdot \sum_{i=1}^d	2|\epsilon^+_{i,j,p(\kappa)}-\epsilon^-_{i,j,p(\kappa)}| \cdot \sup\{ N'(\tfrac{1}{t}K+v) \mid v\in[0,1)^d \} \\
										&+ \Big( \sum_{j=1}^{2^{(n^d)}} |w_j|\Big) \cdot t^{d-1} \cdot \sup\{ N'(\tfrac{1}{t}K+v) \mid v\in[0,1)^d \}\end{aligned}				
\\ 
& \le  \Big( \sum_{j=1}^{2^{(n^d)}} |w_j|\Big) \cdot t^{d-1} \cdot \Big( (2nd\mu+1)\cdot \sup\{N'(\tfrac{1}{t}K+v) \mid v\in[0,1)^d \} \Big)
\\
& \le st														
\end{align*}
for an appropriate constant $s$ by Lemma \ref{l:N'bounded} or Lemma \ref{T:N'bounded}. The second assertion is an immediate consequence of the first one. \qed\medskip



\section{Counterexamples} \label{sec:counter}

In this section we give examples showing that we cannot relax the assumptions of Theorem \ref{T:main}. At first we show that without the assumption that the different boundary parts $F_1, \dots, F_m$ in (R\ref{item:R_R2}) do not intersect under an angle of zero the convergence rate of $\Var \hat S_t(K+tU)$ may be slower. Building Cartesian products with $[0,1]^{d-2}$ one obtains a similar result in higher dimensions. It remains open, whether convergence can be assured at all without this assumption.

\begin{exa}\label{Ex:order}
Put for some $k\ge 2$ (see Figure \ref{F:order})
\begin{align*} 
L_1 &:= \{ (x,y) \in\mathbb{R}^2 \mid  |x|^k \le y \le 1 \} \\
L_2 &:= \{ (x,y) \in\mathbb{R}^2 \mid -1\le x \le 1,\, -1 \le y\le 0 \}\mbox{ and} \\
K &:= L_1 \cup L_2.
\end{align*}

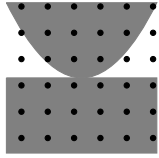
\begin{figure}
\begin{center}
\begin{tikzpicture}

\filldraw[gray] (-1,0)--(-1,-1)--(1,-1)--(1,0)--(-1,0) ;
\filldraw[gray]	(-1,1) parabola bend (0,0) (1,1)--(-1,1);

\foreach \x in {0.95, 0.6, 0.25, -0.1, -0.45, -0.8}
\foreach \y in {0.95, 0.6, 0.25, -0.1, -0.45, -0.8}
\filldraw (\x,\y) circle(1pt);

\end{tikzpicture}
\end{center}

\caption{The set $K$ from Example \ref{Ex:order} in the case $k=2$}
\label{F:order}
\end{figure}

The pixel configuration counts of $K$ can be computed -- except for the complete white pixel configuration -- by
\[ N_{t,j}(K+tU) = N_{t,j}(L_1+tU) + N_{t,j}(L_2+tU) -E_j, \]
where $U=(U_1,U_2)$ is uniformly distributed on $[0,1)^2$ and where $E_j$ is some correction term arising due to intersection effects. Assume $n=2$ from now on, i.e.\ consider $2\times 2$-pixel configurations. For the pixel configuration $(B_j,W_j)$ consisting of two white pixels in the lower horizontal row and two black pixels in the upper horizontal row this correction term is  
\[E_j= D:= \# \{ v_1\in\mathbb{Z} \mid (tv_1,t), (tv_1+t,t) \in L_1+tU \}. \]
For the opposite pixel configuration $E_j=D+2$ has to be subtracted, for the complete black pixel configuration $E_j=-D$ holds, for the two pixel configurations which contain two white pixels in the lower row and one black pixel and one white pixel in the upper row we have $E_j=1$  and for the the two pixel configurations with black lower row and one black pixel and one white pixel in the upper row, $E_j=-1$ holds. For any other pixel configuration we have $E_j=0$. 

By Theorem \ref{T:main} both $\Var( t N_{t,j}(L_1+tU))$ and $\Var( t N_{t,j}(L_2+tU))$ are of order $O(t^2)$ as $t\to 0$. 

We have
\begin{align*} 
D &= \# \{ v_1\in\mathbb{Z} \mid |tv_1-tU_1|^k \le t-tU_2, |tv_1+t-tU_1|^k \le t-tU_2 \} \\
&=\# \{v_1\in\mathbb{Z} \mid |tv_1-tU_1|\le \sqrt[k]{t-tU_2}, \, |tv_1+t-tU_1| \le \sqrt[k]{t-tU_2} \}
 \begin{cases} \le \lfloor 2\cdot \sqrt[k]{t-tU_2}/t \rfloor \\  \ge \lfloor 2\cdot \sqrt[k]{t-tU_2}/t \rfloor -1.\end{cases}
\end{align*}
Hence
\begin{align*}
\mathbb{P}(D >\rho^*) \ge \int_0^1 \mathbf{1}_{(\rho^*,\infty)}\big(\lfloor 2 \sqrt[k]{t-tu}/t\rfloor-1 \big) \, du \ge \int_0^1 \mathbf{1}_{(0,\frac{1}{3})}(u) \, du =\frac{1}{3}
\end{align*}
for $\rho^*= 2\sqrt[k]{\frac{2}{3}t}/t - 2$ and, similarly,
\[\mathbb{P}( D < \rho_*) \ge \frac{1}{3} \]
for $\rho_* = 2\sqrt[k]{\frac{1}{3}t}/t$.

With $\mu:= \mathbb{E} D$ we get
\[ \frac{1}{3} \le \mathbb{P}\Big( |D-\mu| >\frac{\rho^*-\rho_*}{2} \Big) \le \frac{4\Var D}{(\rho^*-\rho_*)^2} = \frac{ \Var D}{\Big( \Big(\sqrt[k]{\tfrac{2}{3}}-\sqrt[k]{\tfrac{1}{3}} \Big) \cdot t^{-1+\frac{1}{k}} -1 \Big)^2} \]
and therefore $\Var D\ge c_k \cdot t^{-2+2/k}$ for sufficiently small $t>0$ for an appropriate constant $c_k>0$. 
This yields that $\Var N_{t,j}(K+tU)$ is (at least) of order $t^{-2+\frac{2}{k}}$ for the two pixel configurations $(B_j,W_j)$ that consists of one black row and one white row. The resulting estimator has a variance of order $t^{2/k}$ (unless the weights of these two pixel configurations sum up to zero or \eqref{e:weights_zero} is violated, which is not the case for any reasonable surface area estimator). 
\end{exa}

It is not clear whether Theorem \ref{T:main} still holds, when the regularity assumption that the convex sets $L_1, \dots, L_m$ in (R\ref{item:R_Rd}) are either polytopes or have a $C^2$ boundary with nowhere vanishing Gauss-Kronecker curvature is removed. However, our method of proof breaks down without this assumption. In fact, the following example shows that Lemma \ref{T:N'bounded} does not hold for an arbitrary compact and convex set $K$.

\begin{exa}\label{Ex:polytop}
Let $K\subseteq \mathbb{R}^3$ be (see Figure \ref{F:polytop}) the convex hull of
\[ \big\{ (\sin\big(\tfrac{1}{2k}\big),\cos\big(\tfrac{1}{2k}\big),1\big) \mid k\in\mathbb{N} \big \} \cup 
	\big\{ (\sin\big(\tfrac{1}{2k-1}\big),\cos\big(\tfrac{1}{2k-1}\big),0\big) \mid k\in\mathbb{N} \big \}. \]
\begin{figure}

\begin{center}
\begin{tikzpicture}

\draw (0,0) circle (0.75*2cm);
\filldraw[lightgray] (0,0)--( 0 , -2 *0.75)--( 2*0.75 , 0 )--( 1.7321*0.75 , 1*0.75 )--( 1.4142*0.75 , 1.4142*0.75 )--( 1.1756*0.75 , 1.618*0.75 )--( 1*0.75 , 1.7321*0.75 )--( 0.8678*0.75 , 1.8019*0.75 )--( 0.7654*0.75 , 1.8478*0.75 )--( 0.684*0.75 , 1.8794*0.75 )--( 0.618*0.75 , 1.9021*0.75 )--(0,2*0.75);
\filldraw[pattern=horizontal lines] ( 0 , 2*0.75 )--( 1.7321*0.75 , -1*0.75 )--( 1.9021*0.75 , 0.618*0.75 )--( 1.5637*0.75 , 1.247*0.75 )--( 1.2856*0.75 , 1.5321*0.75 )--( 1.0813*0.75 , 1.6825*0.75 )--( 0.9294*0.75 , 1.7709*0.75 )--( 0.8135*0.75 , 1.8271*0.75 )--( 0.7225*0.75 , 1.8649*0.75 )--( 0.6494 *0.75, 1.8916*0.75 )--(0,2*0.75);
\foreach \x in {( 1.7321*0.75 , -1*0.75 ), ( 1.9021 *0.75, 0.618*0.75 ), ( 1.5637*0.75 , 1.247*0.75 ), ( 1.2856 *0.75, 1.5321 *0.75), ( 1.0813 *0.75, 1.6825 *0.75), ( 0.9294*0.75 , 1.7709 *0.75), ( 0.8135*0.75 , 1.8271*0.75 ), ( 0.7225*0.75 , 1.8649*0.75 ), ( 0.6494*0.75 , 1.8916*0.75 ), (0.5896*0.75, 1.9112*0.75), (0.5396*0.75, 1.9258*0.75), (0.4974*0.75, 1.9372*0.75), (0.4612*0.75, 1.9460*0.75), (0.4300*0.75, 1.9532*0.75), (0.4026*0.75, 1.9590*0.75), (0.3786*0.75, 1.9638*0.75), (0.3572*0.75, 1.9678*0.75), (0.3380*0.75, 1.9712*0.75), (0.3208*0.75, 1.9742*0.75), ( 0.3052*0.75 , 1.9766*0.75 ), ( 0.2458*0.75 , 1.9848*0.75 ), ( 0.2056*0.75 , 1.9894 *0.75), ( 0.1768 *0.75, 1.9922*0.75 ), ( 0.155*0.75 , 1.994 *0.75), ( 0.138 *0.75, 1.9952*0.75 ), ( 0.1244*0.75 , 1.9962 *0.75), (0,2*0.75)
}
\draw[thick] \x circle(1mm);
\foreach \y in {( 0 , -2*0.75 ), ( 2 *0.75, 0 ), ( 1.7321*0.75 , 1*0.75 ), ( 1.4142 *0.75, 1.4142*0.75 ), ( 1.1756 *0.75, 1.618*0.75 ), ( 1*0.75 , 1.7321*0.75 ), ( 0.8678*0.75 , 1.8019*0.75 ), ( 0.7654*0.75 , 1.8478*0.75 ), ( 0.684*0.75 , 1.8794*0.75 ), ( 0.618*0.75 , 1.9021*0.75 ), ( 0.5176*0.75 , 1.9318*0.75 ), ( 0.445*0.75 , 1.9498*0.75 ), ( 0.3902*0.75 , 1.9616*0.75 ), ( 0.3472*0.75 , 1.9696*0.75 ), ( 0.3128*0.75 , 1.9754*0.75 ), ( 0.3128*0.75 , 1.9754*0.75 ), ( 0.2506*0.75 , 1.9842*0.75 ), ( 0.209*0.75 , 1.989*0.75 ), ( 0.1792*0.75 , 1.992*0.75 ), ( 0.157*0.75 , 1.9938*0.75 ), ( 0.1396*0.75 , 1.9952*0.75 ), ( 0.1256*0.75 , 1.996*0.75 ), (0,2*0.75)}
\filldraw[gray] \y circle(0.8mm);
\end{tikzpicture}
\end{center}

The intersection of the set $K$ with the level $z=1$ is plotted gray, while the intersection of $K$ with the level $z=0$ is displayed shaded. The vertices at level $z=1$ are marked by solid gray dots, while the vertices at level $z=0$ are marked with empty black dots. 

\caption{The set $K$ of Example \ref{Ex:polytop}}
\label{F:polytop}
\end{figure}

Consider the homothetic image $rK$ of $K$. It has edges $l(r, k)$, $k\in\mathbb{N},$ connecting $\big(r\sin\big(\tfrac{1}{2k}\big),r\cos\big(\tfrac{1}{2k}\big),r\big)$  to $\big(r\sin\big(\tfrac{1}{2k-1}\big),r\cos\big(\tfrac{1}{2k-1}\big),0\big)$. Clearly each edge $l(r,k)$ separates two different boundary components $H_p$ and $H_{p'}$, $p\ne p'$. Thus a cell whose interior intersects an edge $l(r,k)$ belongs to $\mathcal{C}'(rK)$. Moreover, the numbers of cells necessary to cover one edge $l(r,k)$ grows asymptotically linear in $r$. If $r$ is large enough depending on $k$, then a cell intersecting one edge of the form $l(r,k)$ cannot intersect any other edge $l(r,k')$, $k'\ne k$. Thus the number of cells necessary to cover all lines $l(r,k), k\in\mathbb{N},$ must grow faster than linear. We conclude that $N'(rK)$ grows faster than linear.

\end{exa}

\section{Simulations} \label{sec:sim}

In this section we evaluate the variances of local estimators based on simulations. We use the weights obtained from discretizing the Crofton formula; see \cite{SON06} or \cite[Sec. 5.2]{OS09}.

\begin{figure}%
\includegraphics[width=0.45\columnwidth]{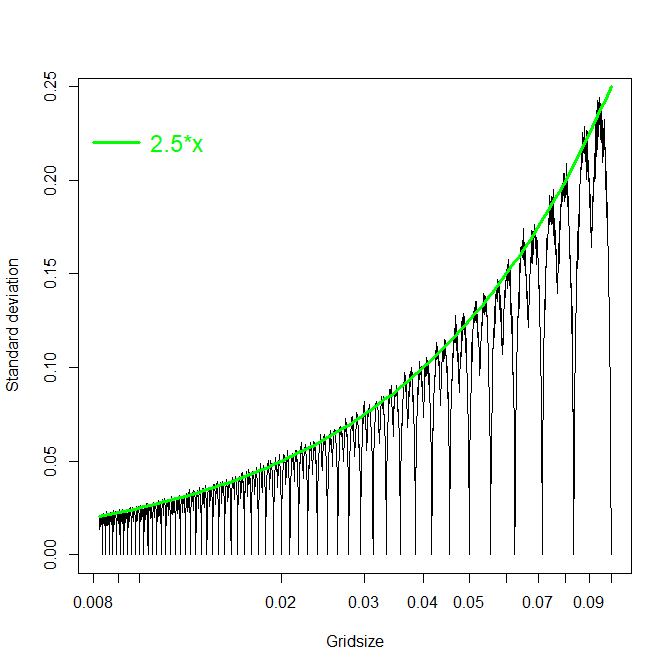}%
\includegraphics[width=0.45\columnwidth]{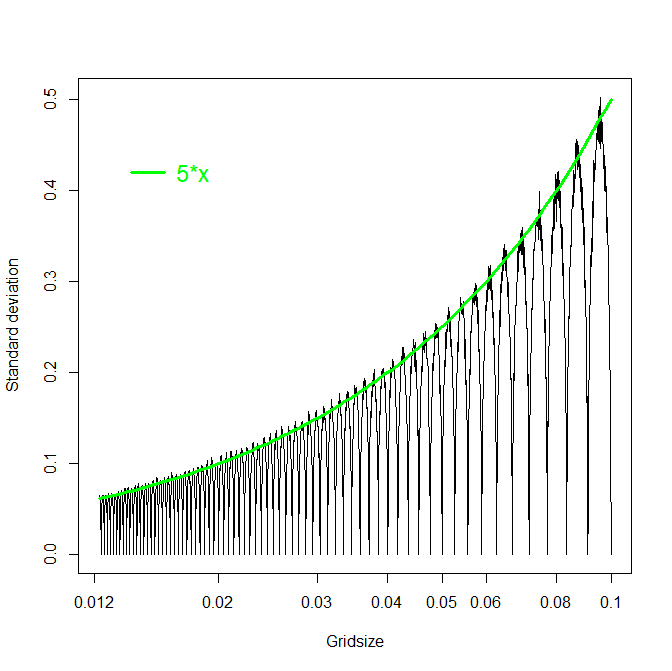}

\includegraphics[width=0.45\columnwidth]{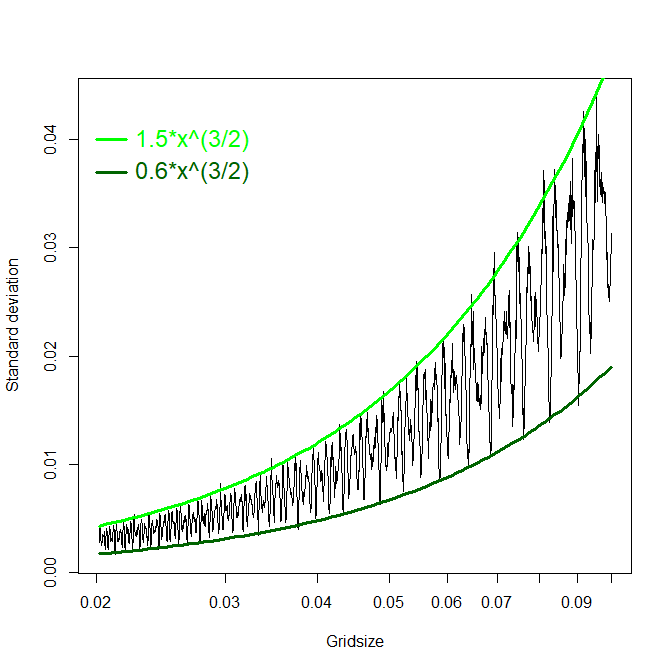}

The standard deviation of $\hat S_t(K+tU)$ in dependence of the grid size $t$. In the upper left picture $K$ is the cuboid mentioned at the beginning of Section \ref{sec:sim}, in the upper right picture $K$ is the skewed parallelepiped and in the lower picture $K$ is the ball. For further details see the text. 

\caption{The standard deviation of the surface area estimator}%
\label{F:sim_surface}%
\end{figure}

We consider three different objects: A cuboid with axes parallel to the coordinate axes and side-lengths $\frac{1}{2}$, $1$ and $1$, a parallelepiped with vertices $(0, 0, 0)$, $(1, -1, 0)$, $(-1, 2, 0)$, $(0, 1, 0)$, $(0, 0, 1)$, $(1, -1, 1)$, $(-1, 2, 1)$ and $(0, 1, 1)$ and a ball of radius $1$. We evaluate the variance at lattice distances of the form $0.1 \cdot 0.999^k$ for integers $k$. However, we replace for each integer $l$ the largest grid size which is smaller than $1/l$ by $1/l$, since we expect the variance to have  local minimum at these points for reasons explained below. For each object at each grid size we determine the variance of the surface area estimator using $400$ simulation runs. The results are reported in Figure \ref{F:sim_surface}.

For both parallelepipeds we see a highly oscillatory behavior of the standard deviation of $\hat S_t(K+tU)$. It drops down to zero if $t=2/k$ resp.\ $t=1/k$. This is explained by the fact that for such $t$ the distance of two opposite sides of the parallelepiped is an integer multiple of the grid size, hence a point of the pixel lattice enters the set $K+tU$ exactly at the time another point leaves at the opposite side (we imagine that $U$ is varying) and thus the pixel configuration counts $N_{t,j}(K+tU)$ do not depend on $U$. We see that the upper bound is a linear function -- as we would expect from Theorem \ref{T:main} (in the pictures it does not look linearly, since the $x$-axis is logarithmically scaled).   

For the ball we do not see regular oscillations anymore but instead we see a Zitterbewegung. This is not surprising, since a similar behavior of the variance is known for volume estimates based on pixel counts, see e.g.\ \cite{Ma71,Mat89}. As already observed by Lindblad \cite{Li05}, both the upper and the lower bound behave approximately as some constant times $t^{3/2}$. This indicates that Theorem \ref{T:main} does not provide the best possible bound in the case that $K$ is a ball.     

\begin{figure}%
\includegraphics[width=0.45\columnwidth]{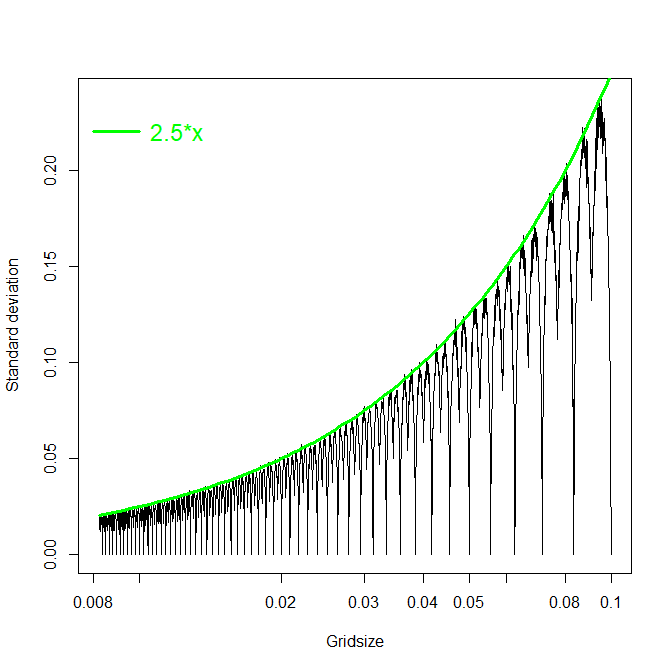}%
\includegraphics[width=0.45\columnwidth]{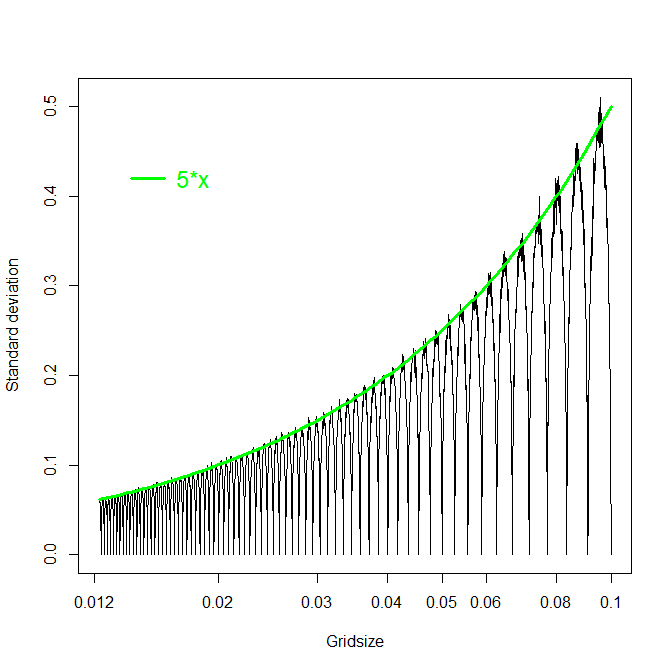}

\includegraphics[width=0.45\columnwidth]{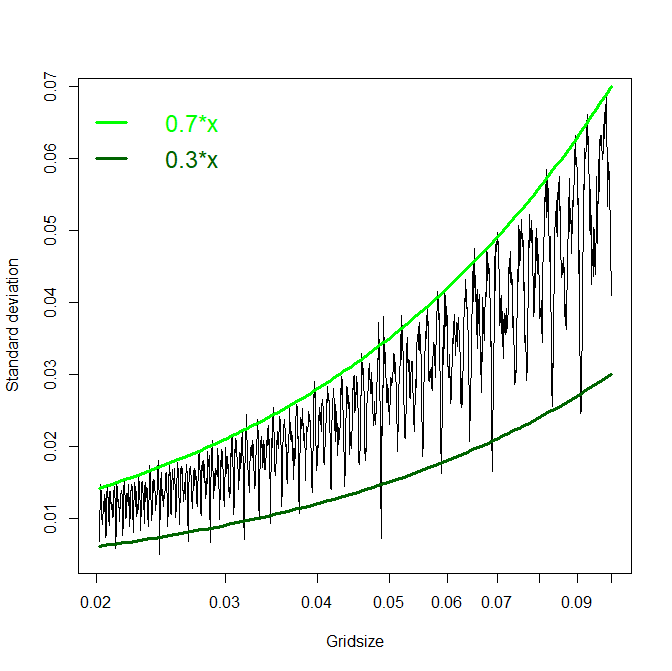}

The standard deviation of the estimator of the integral of mean curvature from \cite{SON06} applied to $K+tU$ in dependence of the grid size $t$. In the upper left picture $K$ is the cuboid mentioned at the beginning of Section \ref{sec:sim}, in the upper right picture $K$ is the skewed parallelepiped and in the lower picture $K$ is the ball. For further details see the text. 

\caption{The standard deviation of the estimator of the integral of mean curvature}%
\label{F:sim_intmecu}%
\end{figure}

The surface area is (up to a factor $\frac{1}{2}$) the $(d-1)$-st intrinsic volumes. The intrinsic volumes (or Minkowski functionals) on $\mathbb{R}^d$ are a family of $d+1$ geometric functionals, including beside the surface area also the volume, the integral of mean curvature and the Euler characteristic. In $\mathbb{R}^2$ and $\mathbb{R}^3$ there are no further intrinsic volumes and in $\mathbb{R}^2$ the surface area and the integral of mean curvature coincide. The variances of volume estimators have been studied intensively, see \cite{IKKN06} for an overview and \cite{Gu15} for a more recent development. The Euler characteristic is a purely topological quantity. A sufficiently smooth set $K$ can be reconstructed up to homeomorphism from a pixel image of sufficiently fine resolution \cite{Pa82, SLS07} and thus it is possible to construct an estimator for the Euler characteristic that returns the correct value with probability one. A natural question is how estimators for the remaining intrinsic volume in $\mathbb{R}^3$, i.e.\ the integral of mean curvature, behave. We have simulated the standard deviation of the estimator for the integral of mean curvature from \cite{SON06, OS09} under the same setup as above. The results are shown in Figure \ref{F:sim_intmecu}. 

For the parallelepipeds  the plots in Figure \ref{F:sim_intmecu} look quite similar to the plots in Figure \ref{F:sim_surface}. This is partially explained by the fact that the breakdowns to zero of the standard deviation are caused by the same effects and therefore take place at the same lattice distances. However, one sees that even the scales are the same. This is definitely coincidence -- if one replaces $K$ by a homothetic image $\lambda K$ the scales change differently for the estimator of the integral of mean curvature than for the surface area estimator. Again, one sees a linear decrease of the maximal standard deviation as the lattice distance tends to zero -- however, the maximal standard deviation now decreases only linearly for the sphere as well.

\section{Discussion and open questions} \label{sec:discuss}

We have shown that the variances of local estimators for the surface area from binary pixel images are of order $O(t^2)$ as the lattice distance $t$ tends to zero. So they are asymptotically neglectable compared to the biases. 

While the situation somewhat reminds of a dart player always hitting approximately the same point which is however far from the middle of the target, this is essentially good news: In $\mathbb{R}^3$ it was clear from the results of \cite{KiZi10} that we cannot derive an upper bound of less than $0.04^2$ for the relative asymptotic mean squared error of any local surface area estimator. Now we have shown that this lowest imaginable upper bound is indeed correct for some estimator.   

It is instructional to review what we have done in the light of the article of Kiderlen and Rataj \cite{KiRa06}. 

\begin{rem}
From Corollary \ref{K:pixel_bounds} we conclude with help of Lemma \ref{l:N'bounded} and Lemma \ref{T:N'bounded} that
\[ \lim_{t\to 0} \mathbb{E}\, \hat S_t(K+tU) = \sum_{j=1}^{2^{(n^d)}} w_j \sum_{\kappa=1}^\mu \sum_{i=1}^d (\epsilon^+_{i,j,p(\kappa)} - \epsilon^-_{i,j,p(\kappa)}) \cdot I_{\kappa,i}. \]
According to \cite[Theorem~5]{KiRa06} this limit equals (in the case that $U$ is distributed uniformly on $[0,1]^d$) 
\[ \sum_{j=1}^{2^{(n^d)}} w_j \int_{S^{d-1}} (- h(B_j\oplus \check W_j, u))^+ S_{d-1}(K,u), \]
where $S_{d-1}(K,\cdot)$ is the surface area measure of $K$ and $\check W:=\{-w \mid w \in W \}$. 
This shows that the two expressions are equal, but of course this is not the most efficient way of deriving this equality. 
\end{rem}

As an alternative to the estimators examined in this paper, one can also use the surface area estimators based on gray-scale images proposed by Svane \cite{Sv14b}. The biases of these estimators converge to zero \cite{Sv14b}, while the biases of the local estimators for binary images converge to positive values. The variances of the estimators from \cite{Sv14b} converge to zero of order $O(t^{d-1})$ as $t\to 0$ \cite{Sv15}, while the variance of the estimators studied in the present paper converge to zero of order $O(t^2)$. So in $\mathbb{R}^2$ we have a bias-variance trade-off, while in $\mathbb{R}^3$ the variances of both kinds of estimators have the same order of convergence. Moreover, according to the mean squared error Svane's estimators perform better in any dimension. However, these estimators rely on quite severe assumptions. First the image $A$ is assumed to be the convolution of the displayed object $K$ with the point spread function. This is not a reasonable assumption if the image is recorded using computed tomography. Second the point spread function is assumed to be known. So the overall recommendation for practitioners is the following: If the severe assumptions on which the estimators from \cite{Sv14b} rely are fulfilled, then use one of them. However, if these assumptions are not fulfilled, then binarize the image and use an estimator for binary images. 

An interesting question is how the variance of a local estimator of the surface area behaves under other model assumptions. When additionally to the random shift a random rotation is applied to the observed set, the relative mean squared error stays bounded, but does not converge to zero, as the grid distance tends to zero \cite{Li05}. At least for asymptotically unbiased local estimators -- under these model assumptions they exist -- the variance behaves the same way. The question how the variance of a local estimator behaves when it is applied to a Boolean model is open.

The results of Section \ref{sec:sim} show that the asymptotic variances of estimators for the integral of mean curvature in $\mathbb{R}^3$ tend to zero of order $O(t^2)$ as $t\to 0$ as well. However, a proof following the proof of Theorem \ref{T:main} shows only that this variance stays bounded. The reason is that beside \eqref{e:weights_zero} there are several other equalities that weights for any reasonable estimator of the integral of mean curvature have to fulfill. Clearly, these relations would have to be exploited in a proof of the optimal bound. However, up to now they are only known in $\mathbb{R}^2$; see \cite{Sv14z}. Of course, a theoretical determination of the optimal $O$-class of $\Var \hat S_t(K+tU)$ in the special case that $K$ is a ball is desirable as well. Also the variance of the surface area estimator applied to bodies with boundary components that intersect under an angle of zero should be further investigated.

The result that the asymptotic variance of an estimator $\hat S_t$ of the form \eqref{e:form_local} is asymptotically neglectable compared to its bias can also be interpreted as $\hat S_t$ being in fact not an estimator for the surface area, but for some other geometric quantity.  For compact, convex sets $K\subseteq \mathbb{R}^d$ with interior points this quantity is the difference of mixed volumes
\[ V(\mathcal{L}_1[1], K[d-1]) - V(\mathcal{L}_2[1], K[d-1]) \]
for
\[ \mathcal{L}_1 := \sum_{j=1}^{2^{n^d}} w_j\conv\big(\{0\} \cup (B_j\oplus \check W_j) \big) \]
and
\[ \mathcal{L}_2 := \sum_{j=1}^{2^{n^d}} w_j\conv (B_j\oplus \check W_j), \]
where $\sum$ denotes the Minkowski sum and $\conv S$ denotes the convex hull of $S\subseteq \mathbb{R}^d$. 
A Miles-type formula for this quantity has been obtained in \cite{ONS09}. Moreover, it also satisfies the assumptions of \cite{HM99} and thus we have central limit theorems for it applied to germ-grain models. We expect that also the local estimators for the other intrinsic volumes are in fact estimators for new geometric quantities.


\begin{thebibliography}{99}
\bibitem{CFTT03} D.~Coeurjolly, F.~Flin, O.~Teytaud and L.~Tougne: \emph{Multigrid convergence and surface area estimation}, in T.~Asano et al. (eds.): Geometry, Morphology and Computational Imaging, Springer (2003), 101--119.  
\bibitem{Gu15} J.~Guo: \emph{Lattice points in rotated convex domains}, Revista Matem\'atica Iberoamericana {\bf 31} (2015), 411--438.
\bibitem{IKKN06} A.~Ivi\'c, E.~Kr\"atzel, M.~K\"uhleitner and W.~Nowak: \emph{Lattice points in large regions and related arithmetic functions: Recent developements in a very classical topic}, in W.~Schwarz et al. (eds.): Elementare und analytische Zahlentheorie - Proceeding of the 3rd Conference, Franz Steiner Verlag Stuttgart (2006), 89--128.    
\bibitem{JK10} J.~Jan\'a\u{c}ek and L.~Kub\'inov\'a: \emph{Variances of length and surface area estimates by spatial grids: preliminar study}, Image Analysis \& Stereology {\bf 29} (2010), 45--52.  
\bibitem{HS89} U.~Hahn and K.~Sandau: \emph{Precision of surface area estimation using spatial grids}, Acta Stereologica {\bf 8} (1989), 425--430. 
\bibitem{HM99} L.~Heinrich and I.~Molchanov: \emph{Central limit theorem for a class of random measures associated with germ-grain models}, Advances in Applied Probability (SGSA) {\bf 31} (1999), 283--314.
\bibitem{KiRa06} M.~Kiderlen and J.~Rataj: \emph{On infinitesimal increase of volumes of morphological transforms}, Mathematika {\bf 53} (2006), 103--127. 
\bibitem{KlRo04} R.~Klette and A.~Rosenfeld: \emph{Digital Geometry}, Elsevier (2004). 
\bibitem{KlSu01} R.~Klette and H.~Sun: \emph{Digital planar segment based polyhedrization for surface area estimation}, in C.~Arcelli et al.\ (eds.): 4th International Workshop on Visual Form (2001), 356--366.
\bibitem{Li05} J.~Lindblad: \emph{Surface area estimation of digitized 3d objects using weighted local configurations}, Image and Vision Computing {\bf 23} (2005), 111--122. 
\bibitem{LiNy02} J.~Lindblad and I.~Nystr\"om: \emph{Surface area estimation of digitized 3d objects using local computations}, in A.~Braquelaire et al.\ (eds.): 10th International Conference on Discrete Geometry for Computer Imagery (2002), 267--278.
\bibitem{Mat89} B.~Mat\'ern: \emph{Precision of area estimation: a numerical study}, Journal of Microscopy {\bf 153} (1989), 269--284. 
\bibitem{Ma71} G.~Matheron: \emph{The Theory of Regionalized Variables and its Applications}, Les Cahiers du Centre de Morphologie Math\'ematique de Fontainebleau (1971).  
\bibitem{ONS09} J.~Ohser, W.~Nagel and K.~Schladitz: \emph{Miles formulae for Boolean models observed on lattices}, Image Analysis \& Stereology {\bf 28} (2009), 77--92.
\bibitem{OS09} J.~Ohser and K.~Schladitz: \emph{3d Images of Material Structures}, Wiley, Weinheim (2009). 
\bibitem{Pa82} T.~Pavlidis: \emph{Algorithms for Graphics and Image Processing}, Computer Science Press (1982). 
\bibitem{SON06} K.~Schladitz, J.~Ohser and W.~Nagel: \emph{Measuring intrinsic volumes in digital 3d images}, in A.~Kuba et.\ al.\ (eds.): 13th International Conference on Discrete Geometry for Computer Imagery (2006), 247--258.
\bibitem{Schn14} R.~Schneider: \emph{Convex Bodies - The Brunn-Minkowski Theory}, Cambridge University Press (2014). 
\bibitem{SW08} R.~Schneider and W.~Weil: \emph{Stochastic and Integral Geometry}, Springer (2008). 
\bibitem{SLS07} P.~Stelldinger, L.~Latecki and M.~Siqueira: \emph{Topological equivalence between a 3d object and the reconstruction of its digital image}, IEEE Transactions on Pattern Analysis and Machine Intelligence {\bf 29} (2007), 126--140.
\bibitem{Sv14z} A.~Svane: \emph{Local digital estimators of intrinsic volumes for Boolean models and in the design-based setting}, Advances in Applied Probability (SGSA) {\bf 46} (2014), 35--58.
\bibitem{Sv14} A.~Svane: \emph{On multigrid convergence of local algorithms for intrinsic volumes}, Journal of Mathematical Imaging and Vision {\bf 49} (2014), 148--172.
\bibitem{Sv14b} A.~Svane: {\it Estimation of intrinsic volumes from digital grey-scale images}, Journal of Mathematical Imaging and Vision {\bf 49} (2014), 352--376.
\bibitem{Sv15} A.~Svane: {\it Asymptotic variance of grey-scale surface area estimators}, Advances in Applied Mathematics {\bf 62} (2015), 41--73.
\bibitem{KiZi10} J.~Ziegel and M.~Kiderlen: {\it Estimation of surface area and surface area measure of three-dimensional sets from digitizations}, Image Vision and Computing {\bf 28} (2010), 64--77.
\end{thebibliography}
\end{document}